\documentclass{article}
\usepackage[%
journal=MSL,    %% Replace XXX by one of the above journal codes.
lang=british,   %% Change to `american' if you use American English.
]{ems-journal}

\usepackage{amsfonts}
\usepackage{algorithmic}
\usepackage{algorithm}
\usepackage{textcomp}
\usepackage{stfloats}
\usepackage{verbatim}

\numberwithin{equation}{section}

\theoremstyle{Lemma}
\newtheorem{lemma}{Lemma}

\theoremstyle{Corollary}
\newtheorem{corollary}{Corollary}

\theoremstyle{Theorem}
\newtheorem{theorem}{Theorem}

\theoremstyle{Remark}

\begin{document}

%------
% Insert the title of your paper and (if necessary)
% a short title for the running head.
%------
\title{Rapid Bayesian Computation and Estimation for Neural Networks via Log-Concave Coupling}
\titlemark{Rapid Bayesian Computation and Estimation for Neural Networks via Log-Concave Coupling}

%------

%%%% Pls fill in all fields for each author
%%%% Label the authors by their position in the authors' list using {}
%%%% If you published any math paper ever, you have an MR Author ID.
%  Please look it up in three easy (and free) steps:
% 1. copy the bibliographic data of any published paper (co-)authored by you in the search field at https://mathscinet.ams.org/mathscinet/freetools/mref
% 2. Hit your name in the search result
% 3. Find your MR Author ID in the first row, copy it in the \mrid{} field
%%%% If you have not created your ORCID yet, you may like to do it now, pls copy it in the field \orcid{}
%%%% Abbreviate first names for the running head

\emsauthor{1}{
	\givenname{Curtis}
	\surname{McDonald}
	\mrid{1398664}
	\orcid{0000-0001-6298-2861}}{C. McDonald}
%%%% Repeat the same fields for each numbered author
\emsauthor{2}{
	\givenname{Andrew R}
	\surname{Barron}
	\mrid{222645}
	\orcid{}}{A. R. Barron}

%%%% Please provide detailed address info for each author
%%%% Use the same numbering as for \emsauthor above
%%%% Please look up the ROR ID of your institute here: https://ror.org
\Emsaffil{1}{
	\department{Simons Institute for the Theory of Computing}
	\organisation{University of California, Berkeley}
	\rorid{03v76x132}
	\zip{94720-2190}
	\city{Berkeley}
	\country{United States of America}
	\affemail{c.mcdonald.simons@berkeley.edu}}
%%%% Repeat the same fields for each numbered author
%%%% If some author has multiple affiliations, repeat the fields for each affiliation
%%%% Number the affiliations using {}
\Emsaffil{2}{
	\department{Department of Statistics and Data Science}
	\organisation{Yale University}
	\rorid{03v76x132}
	\address{P.O. Box 208290}
	\zip{06520-8290}
	\city{New Haven}
	\country{United States of America}
	\affemail{Andrew.Barron@yale.edu}}

%------
% Add MSC 2020 codes according to https://zbmath.org/classification/.
% A unique primary MSC code (in curly brackets) is mandatory,
% while secondary MSC codes (in square brackets) are optional.
%------
\classification[62M45, 65C05]{62F15}

%------
% Add a list of keywords.
%------

%------
% Insert your abstract.
%------

\begin{abstract}
This paper presents the study of a Bayesian estimation procedure for single-hidden-layer neural networks using $\ell_{1}$ controlled neuron weight vectors. We study the structure of the posterior density and provide a representation that makes it amenable to rapid sampling via Markov Chain Monte Carlo (MCMC). Let the neural network have $K$ neurons with internal weights of dimension $d$ and fix the outer weights. Thus there are $Kd$ parameters overall. With $N$ data observations, use a gain parameter or inverse temperature of $\beta$ in the posterior density for the internal weights.

The posterior is intrinsically multi-modal and not naturally suited to rapid mixing of direct MCMC algorithms. For a continuous uniform prior on the $\ell_{1}$ ball, we demonstrate that the posterior density can be written as a mixture density with suitably defined auxiliary random variables, where the mixture components are log-concave. Furthermore, when the total number of model parameters $Kd$ is large enough that $Kd \geq C(\beta N)^{2}$, the mixing distribution of the auxiliary random variables is also log-concave. Thus, neuron parameters can be sampled from the posterior by only sampling log-concave densities. The authors refer to the pairing of weights with such auxiliary random variables as a log-concave coupling.
\end{abstract}

\maketitle

\tableofcontents

%------
\section{Introduction}
Single-hidden-layer artificial neural networks provide a flexible class of parameterized functions for data fitting applications. Specifically, denote a single-hidden-layer neural network as the parameterized function
\begin{equation}
    f_{w}(x) = f(w,x) = \sum_{k=1}^{K}c_{k}\psi(w_{k} \cdot x),\label{nn_def}
\end{equation}
with $K$ neurons, activation function $\psi$, and interior weights $w_{k} \in \mathbb{R}^{d}$. 
Fix a positive scaling $V$ and let the exterior weights $c_{k}$ be positive or negative values $c_{k} \in \{-\frac{V}{K},\frac{V}{K}\}$. Thus, $f_{w}(x)$ is a convex combination of $K$ signed neurons scaled by $V$. Constant and linear terms $c_0 + w_0 \cdot x$ may be added in the definition of $f_w(x)$ to achieve additional flexibility, though we will not address that matter explicitly.

 While the class of functions which can be modeled well by a shallow neuron network is well understood, computational guarantees for algorithms to fit specific neural networks have been lacking. The approximation ability of these networks has been studied for many years, which we briefly review here. Assume input vectors $x \in \mathbb{R}^{d}$ as having bounded entries, $x \in [-1,1]^{d}$, and interior weights have bounded $\ell_{1}$ norm $\|w_{k}\|_{1} \leq 1$. Thus the input to each neuron is between $[-1,1]$. 

 The early work of \cite{barron1993universal} showed that moderately wide single-hidden-layer networks with sigmoid activation functions can accurately approximate target functions with a condition on the Fourier components of the target function. For a sigmoid activation function and $K$ neurons, the squared error with the target function was shown to be on the order of $O(\frac{1}{K})$.
 
Denote the set of signed neurons with $\ell_{1}$ controlled interior weights as the collection of functions $h:[-1,1]^{d} \to \mathbb{R}$
\begin{align}
\Psi = \{h:h(x)=\pm \psi(w \cdot x), \|w\|_{1}\leq 1\}.
\end{align}
The closed convex hull of $\Psi$ includes functions $f$ which can be written as a possibly infinite mixture of signed neurons, and functions which are the limit of a sequence of such mixtures. Specializing the results of \cite{barron1992neural},\cite{barron1993universal},\cite{klusowski2018approximation}, networks of the form \eqref{nn_def} provide accurate approximation for functions $f$ with $\frac{f}{V}$ in the closure of the convex hull of $\Psi$. The infimum of such $V$ is called the variation $V_{f}$ of the function $f$ with respect to the dictionary $\Psi$.

Assume a data set of size $N$ of iid data $(x_{i}, y_{i})_{i=1}^{N}$ and $E[y_{i}|x_{i}] =f(x_i)$ with $f(x)$ in the closure of the convex hull of $\Psi$. Results in \cite{barron2019complexity} show for network of width $K = O([N/\log(d)]^{1/2})$, the empirical risk minimizer (ERM) $f_{\hat{w}}$ yields a statistical risk control of the order
\begin{align}
    E[\|f_{\hat{w}}-f\|^{2}]&=O(\Big(\frac{\log(d)}{N}\Big)^{\frac{1}{2}}),
\end{align}
provided there is sub-Gaussian control of the distribution of the response $Y$. The expectation here is with respect to the training data, while the norm square provides the expectation for the loss at an independent new input vector. Analogous deep net conclusions are also in \cite{barron2018approximation}, \cite{barron2019complexity}.

Given that the ERM can achieve such risk control, there is much interest to understand theoretically the optimization properties of neural networks by gradient-based properties. One line of investigation as in  \cite{cai2024large,tsigler2023benign,frei2024double,jacot2018neural,liu2022loss} compares the network to a certain infinite width limit under initialization and scaling assumptions (called the neural tangent kernel, NTK), where the network trained under gradient descent are shown to approach a kernel regression solution linear in the response vector $y$. They also utilize a scaling of $1/\sqrt{K}$ on their outer weights rather than the $1/K$ scaling we use.

With finite $K<N$ and outer weights that scale as $c_{k} = 1/K$, even in the single-hidden-layer case no known optimization algorithm is able to solve this optimization problem in a polynomial number of iterations in $N$ and $d$. Instead, in this work we move away from an optimization approach to choosing neuron parameters and use a Bayesian method of estimation placing a posterior distribution on neuron parameters.

Bayesian neural networks have been studied for many years \cite{neal_bayesian_1996,gallego2022current,charnock_bayesian_2020}, although specific mixing time bounds for Markov Chain Monte Carlo (MCMC) to guarantee polynomial time complexity have been a barrier to their implementation. Recent approaches have studied the simplification of the posterior in the NTK regime, resulting in the posterior being near the posterior associated with a Gaussian process prior \cite{hron2022wide,hanin_bayesian_2024}. These approaches require $K/N \to \infty$ to achieve that simplification of the posterior density. The bounded $K/N$ setting is shown in \cite{hron2022wide,hanin_bayesian_2024} to be distinct with potentially more flexible non-Gaussian process behavior. Indeed, such flexibility arises in our model where the internal weights are adapted by the posterior.

We quantify when sampling can be accomplished in polynomial time with MCMC. Say we have data consisting of $N$ input and response pairs $(x_{i}, y_{i})_{i=1}^{N}$. Define a prior distribution $P_{0}$ on neuron parameters and a gain or inverse temperature $\beta>0$. Define the posterior density
\begin{align}
    p(w|x^{N}, y^{N})& \propto p_{0}(w)e^{-\frac{\beta}{2}\sum_{i=1}^{N}(y_{i}-f(w, x_{i}))^{2}},
\end{align}
and the associated posterior mean at a given $x$ value  $\mu(x)= E_{P}[f(x,w)|x^{N}, y^{N}]$.
%\begin{align}
%    \mu(x)&= E_{P}[f(x,w)|x^{N}, y^{N}].
%\end{align}

The barrier to implementing the Bayesian approaches defined above is being able to sample from the density  $p(w|x^{N}, y^{N})$ and thus compute the posterior average $\mu(x)$. The densities $p(w|x^{N}, y^{N})$ will be high-dimensional and multi-modal densities with no immediately evident structure that would make sampling possible. 

The neural network model has $Kd$ parameters and we have $N$ data observations. A natural method to compute the required posterior averages would be a MCMC sampling algorithm. A key issue for an MCMC method is whether it provably gives accurate sampling in a low polynomial number of iterations in $K,d,N$. Any exponential dependence on the parameters of the problem is not practically useful.

The most common sufficient condition for proving rapid mixing of MCMC methods is log-concavity of the target density \cite{bakry1985diffusions, bakry2014analysis, lovasz2007geometry, dwivedi2019log}. As such, we want to find a representation of the problem built from log-concave densities, so that we may restrict our computation task to only require sampling from log-concave densities.

We show that with the use of an auxiliary random variable $\xi$, the posterior density $p(w|x^{N}, y^{N})$ can be re-written as a mixture density (also called a measure decomposition in the language of \cite{montanari2024provably})
\begin{align}
p(w|x^{N}, y^{N})&= \int p(w|\xi,x^{N}, y^{N})p(\xi|x^{N}, y^{N})d\xi,
\end{align}
using a reverse conditional density $p(w|\xi) = p(w|\xi, x^{N}, y^{N})$ and an induced marginal density $p(\xi) = p(\xi|x^{N}, y^{N})$. When considering a fixed input and response sequence, we will drop the $x^{N}$ and $y^{N}$ conditioning notation. For a certain choice of priors and relationships between $d,K,N,\beta$, we show the reverse conditional is a log-concave density, and the induced marginal for $\xi$ is also a log-concave density. We call such a joint distribution $p(w, \xi)$ that preserves the target marginal $p(w)$ and has a log-concave marginal distribution $p(\xi)$ and a log-concave conditional distribution $p(w|\xi)$ a \textit{log-concave coupling}. As such, samples for $w$ from the posterior can be produced by merely sampling from log-concave densities: that is, we sample from the density of $\xi$ followed by sampling from the density of $w$ given $\xi$. 

For a continuous uniform prior on the $\ell_{1}$ ball, we demonstrate the mixture is a log-concave coupling when the total number of parameters $Kd$ is large enough such that
\begin{align}
    Kd\geq C (\beta N)^{2},
\end{align}
for a given constant $C$ that depends only on the range of data values and the scaling $V$ of the network.

We presume access to a sampling algorithm able to produce accurate samples from a log-concave density in a number of iterations proportional to a low polynomial power of the number of model parameters \cite{dwivedi2019log,lovasz2007geometry,livingstone2019geometry,kook2024gaussian,kook2024sampling}. We leave the specific choice of this algorithm in our setting (e.g. Metropolis Adjusted Langevin Diffusion (MALA), Hamiltonian Monte Carlo, etc.) as well as the tuning of parameters as a technical study for future work, and treat the sampling algorithm as a black box method available to the user. Then one can sample a value for $\xi$ from its marginal, and a value $w|\xi$ from its reverse conditional, resulting in a true draw from the posterior distribution for $w$.

\section{Notation}\label{note_section}
Here we present the mathematical notation used in the paper.
\begin{itemize}
    \item Capital $P$ refers to a probability distribution, while lowercase $p$ is its probability mass or density function.
    \item $f'(\cdot)$ refers to the derivative of a scalar function $f$.
    \item $\nabla$ is the gradient operator and $\nabla^{2}$ is the Hessian operator, producing a matrix of second derivatives.
    \item $\{1, \hdots, N\}$ is the set of whole numbers between 1 and N.
    \item $[a,b]$ is the interval of real values between $a$ and $b$.
    \item $u \cdot v$ is the Euclidean inner product between two vectors.
    \item $u^{\text{\tiny{T}}}, \textbf{X}^{\text{\tiny{T}}}$ refers to the transpose of a vector or matrix, so quadratic forms of a column vector $u$ with a matrix $\textbf{X}$ will be written as $u^{\text{\tiny{T}}}\textbf{X}u$.
    \item $\|w\|_{p}$ refers to the $\ell_{p}$ norm, $\|w\|_{p}=(\sum_{j} (				w_{j})^{p})^{1/p}$.
    \item The $\ell_{1}$ ball is denoted as $S^{d}_{1} = \{w \in \mathbb{R}^{d}:\|w\|_{1} \leq 1\}$.
    \item The $K$ fold Cartesian product of this set is $(S^{d}_{1})^{K}$.
    \item For variables in a sequence, superscripts indicate the set of variables $X^{N} = (X_{i})_{i=1}^{N}$.
    \item Logarithms in the paper are natural logarithms.
\end{itemize}

\section{Bayesian Model and Result}\label{model_section}
Consider input and response pairs $(x_{i}, y_{i})_{i=1}^{N}$ where the $x_{i}$ input vectors are $d$ dimensional and the response values $y_{i}$ are real valued. Consider the $x_{i}$ as being bounded by 1 in each coordinate, $|x_{i,j}| \leq 1$ for all $i \in \{1,\hdots,N\}, j \in\{1, \hdots,d\}$. Further, assume $x_{i,1}=1$ for all $i \in\{1,\hdots,N\}$ so an intercept term is naturally included in the data definition. Accordingly, this requires $d \geq 2$. Denote $\textbf{X}$ as the $N$ by $d$ data matrix which uses the $x_{i}$ as its rows.

Recall the definition of a single-hidden-layer neural network in equation \eqref{nn_def}. We restrict the class of neuron activation functions we consider to have two bounded derivatives with $\psi(0) = 0, |\psi(z)|\leq a_{0}, |\psi'(z)|\leq a_{1}$ and $|\psi''(z)|\leq a_{2}$ for all $z \in [-1,1]$. We assume $a_{0}, a_{1}, a_{2} \geq 1$. This includes for example the squared ReLU activation function $\psi(z) = a (z_{+})^{2}$ or the scaled tanh activation function $\psi(z) = a \tanh(cz)$.

Fix a positive $V$ and let the exterior weights $c_{k}$ be positive or negative values $c_{k} \in \{-\frac{V}{K},\frac{V}{K}\}$. Thus, $f_{w}(x)$ is a convex combination of $K$ signed neurons scaled by $V$. Note, if $\psi$ is odd symmetric as in the case of the tanh activation function, the $c_{k}$ can be all set to positive $\frac{V}{K}$. For non-symmetric activation functions, we can use twice the variation $\tilde{V} = 2V$ and use twice the number of neurons $\tilde{K} = 2K$. For the first $K$ neurons set $c_{k} = \frac{V}{K}$ and for the second set of $K$ neurons set $c_{k} = - \frac{V}{K}$. Under such a structure using $2K$ neurons, any size $K$ network of variation $V$ with any number of positive or negative signed neurons can be constructed from the wider network by setting $K$ of the neurons to be active and the other $K$ to be inactive with $0$ for their interior weight vector. In either case, we consider the outer weights $c_{k}$ to be fixed values, and it is only necessary to train the interior weights $w_{k}$ of the network.

Define $P_{0}$ as a prior measure on $\mathbb{R}^{Kd}$, with density $p_{0}$ with respect to a reference measure $\eta$ (e.g. Lebesgue or counting measure). For each index $i \in \{1, \hdots, N\}$ define the residual of a neural network as
\begin{align}
    \text{res}_{i}(w)&= y_{i}-\sum_{k=1}^{K}c_{k}\psi(w_{k} \cdot x_{i}).
\end{align}
Define the loss function as half the sum of squares of residuals
\begin{align}
    \ell(w)&= \frac{1}{2}\sum_{i=1}^{N}(\text{res}_{i}(w))^{2}.
\end{align}
For any gain parameter $\beta>0$, we define the posterior density 
%\ (with respect to $\beta$) 
as
\begin{align}
\label{p_n_pdf}
    p(w)&= \frac{p_{0}(w)e^{-\beta \ell(w)}}{\int e^{-\beta \ell(w)}p_{0}(w)\eta(dw)}.
\end{align}
Denote the posterior mean with respect to this density as
\begin{align}
    \mu(x)&= E_{P}[f(w, x)].
\end{align}

\subsection{Choice of Prior}
The prior we will consider is the uniform on the set $(S^{d}_{1})^{K}$. That is, independently each weight vector $w_{k}$ is iid uniform on the set of weight vectors with $\ell_{1}$ norm less than 1. This has the density function
\begin{align}
    p_{0}(w)&= \prod_{k=1}^{K}\Big(1\{\|w_{k}\|_{1} \leq 1\} \frac{1}{\text{Vol}(S^{d}_{1})}\Big),
\end{align}
with respect to Lebesgue measure. Note that for each $k$, the vector of absolute values of the coordinates of $w_{k}$ is uniformly distributed on the simplex, which is also a symmetric Dirichlet distribution in $d+1$ dimensions with the all 1's parameter vector.

The authors would like to provide some backstory of why this prior is chosen. The choice of the continuous uniform prior is influenced by previous work using a discrete uniform prior. Consider a discrete version of this density. For some positive integer $M\leq d$, consider the discrete set which is the intersection of $S^{d}_{1}$ with the lattice of points of equal spacing $\frac{1}{M}$. Define this set as $S^{d}_{1, M}$,
\begin{align}
    S^{d}_{1,M}= \{w: M w \in \{-M,\ldots, M\}^{d},~   \|w\|_{1} \leq 1\},
\end{align}
That is, each coordinate $w_{k,j}$ can only be an integer multiple of the grid size $\frac{1}{M}$ and we force the $\ell_{1}$ norm to be less than or equal to 1. We consider the prior under which $w_{k}$ is independent uniform on the discrete set $S^{d}_{1,M}$. This has probability mass function
\begin{align}
    p_{0}(w)&= \prod_{k=1}^{K}\Big(1\{w_{k} \in S^{d}_{1,M}\} \frac{1}{|S^{d}_{1,M}|}\Big),
\end{align}
with respect to counting measure in $(S^{d}_{1,M})^{K}$.

In the author McDonald's PhD thesis \cite{mcdonald2025computation_thesis} various mean-square statistical risk and arbitrary-sequence regret bounds are established for sampling from the posterior distribution with the discrete uniform prior. In particular,  with $M,K = O(N/\log(d))^{1/4})$ and $\beta = (\log(d)/N)^{1/4}$, statistical risk control of the order $O(\log(d)/N)^{1/4}$ is shown with this $\beta$, and faster rates are established with Gaussian distributed $y_i$ and constant $\beta$ equal to the reciprocal of the variance. 
%\ by sampling from the posterior distribution defined by the discrete uniform prior and squared loss with inverse temperature $\beta$. 
 Somewhat larger $M$ and $K$ are allowed, but the risk analysis does not succeed unless $MK$ is of smaller order than $N$. Results of this type are also summarized in \cite{barron2024shannon}, \cite{barron2024log}. 

In this work, we show for the continuous uniform prior there is a log-concave coupling and thus MCMC can be accomplished in polynomial time when $Kd> C(\beta N)^{2}$. With sufficiently large $M$ of larger order than $N$, these mixing time results can be shown for the discrete prior as well. However, this analysis requires $M$ that it is incompatible with the statistical control results. We were not able to combine statistical risk control and polynomial-time sampling simultaneously. In the present work, we only present the polynomial-time sampling results for the continuous prior. 

\subsection{Log-Concave Coupling}
The log-concave coupling result is as follows:

\begin{theorem}\label{paper_main_thm}
Let the neural network have inner weight dimension $d \geq 2$ and $K\geq 2$ neurons with $N$ data observations $(x_{i}, y_{i})_{i=1}^{N}$. Assume $\beta N \geq 2$. Define the values
\begin{align}
    C_{N}&= \max_{i \in\{1, N\}} |y_{i}|+a_{0}V\\
 A_{1}&=2a_{1}+8a_{2}\\
    A_{2}&=2 \sqrt{a_{2}}\\
     A_{3}&=  8\sqrt{e}a_{2}(C_{N}V)^{\frac{3}{2}}[A_{1}+A_{2}(C_{N}V)^{\frac{1}{2}}].
\end{align}
Set a positive $\delta \le 0.05$ and let $d$ and $K$ satisfy,
\begin{align}
K[\log(2Kd)+\log(1/\delta)] &\leq \beta N \label{K_ineq}
\end{align}\par\vspace{-0.5cm}
\noindent and
\par\vspace{-0.8cm}
\begin{align}
    Kd &\geq A_{3}(\beta N)^{2}.\label{strict_ineq}
\end{align}
Using a continuous uniform prior on $(S^{d}_{1})^{K}$ the posterior distribution $p(w)$ can be written as a mixture distribution with an auxiliary random variable $\xi$,
\begin{align}
    p(w)&= \int p(w|\xi)p(\xi)d\xi,
\end{align}
where $p(w|\xi)$ is a log-concave density for each $\xi$, and $p(\xi)$ is a log-concave density. If equation \eqref{strict_ineq} is a strict inequality, $p(\xi)$ is strictly log-concave.
\end{theorem}
\section{Posterior Densities and Log-Concave Coupling}\label{log_concave_coupling}
\subsection{Posterior Density}
Consider the log of the posterior density $p(w)$ defined in equation \eqref{p_n_pdf}, with the continuous uniform prior on $(S^{d}_{1})^{K}$. The log and score of the posterior within the constrained set are equal to the following (with some constant $B$ that is just the log of the normalizing constant)
\begin{align}
    \log p(w)&= -\beta \ell(w)+B\\
    \nabla_{w_{k}} \log p(w)&= \beta\sum_{i=1}^{N}\text{res}_{i}(w)(c_{k}\psi'(w_{k}\cdot x_{i})x_{i}).
\end{align}
Denote the Hessian as $H(w)\equiv \nabla^{2} \log p(w)$. The density $p(w)$ is log-concave if $H(w)$ is negative semi-definite for all choices of $w$. For any vector $u \in \mathbb{R}^{Kd}$, with blocks $u_{k} \in \mathbb{R}^{d}$, the quadratic form $u^{\text{\tiny{T}}}H(w)u$ can be expressed as
\begin{align}
    &-\beta\sum_{i=1}^{N}\Big(\sum_{k=1}^{K} c_{k}\psi'(w_{k}\cdot x_{i})u_{k}\cdot x_{i}\Big)^{2} \label{H_term_1}\\
    &+\beta\sum_{i=1}^{N}\text{res}_{i}(w)\sum_{k=1}^{K}c_{k}\psi''(w_{k} \cdot x_{i})(u_{k}\cdot x_{i})^{2}\label{H_term_2}.
\end{align}
 It is clear that for any vector $u$  the first line \eqref{H_term_1} is a negative term, but term \eqref{H_term_2} may be positive. The scalar values $c_{k}\psi''(w_{k} \cdot x_{i})$ could be either a positive or negative value for each $k$ and $i$, while the residuals $\text{res}_{i}(w)$ can also be positive or negative signed. Thus, the Hessian is not a negative definite matrix in general and $p(w)$ may not be a log-concave density. 

The term \eqref{H_term_2} is capturing how the non-linearity of $\psi$, which provides the benefit of neural networks over linear regression, is complicating matters. If $\psi$ were linear, $\psi''(z)=0$ for all $z$ and we would have a simple linear regression problem. However, since $\psi$ has second derivative contributions, this term must be addressed.

For each data index $i \in \{1, \hdots, N\}$ and each neuron index $ k \in \{1, \hdots,K\}$ we introduce a coupling with an auxiliary random variable $\xi_{i,k}$. The goal of this auxiliary random variable is to force the corresponding individual $i,k$ terms in \eqref{H_term_2} to be negative. Define the values 
\begin{align}
C_{N} &= \max_{i \leq N} |y_{i}|+a_{0}V\\
\rho_{_0,N,K}&= a_{2}\frac{\beta C_{N} V}{K}.\label{rho_def_no_delta}
\end{align}
We denote this $\rho_{0,N,K}$ as $\rho_0$ when it is clear $N$ and $K$ are set.  The role of $\rho_0$ is that of a lower bound on the values of an inverse variance parameter $\rho$ for the conditional distribution of auxiliary random variables below.

Ultimately we will use bounded auxiliary random variables to yield the desired log-concave coupling.  But to motivate the construction, first consider tentatively a simpler unbounded construction.

Conditioning on a weight vector $w$, define the forward coupling as conditionally independent random variables $\xi_{i,k}$ which are normal with mean $w_{k}\!\cdot\!x_{i}$ and variance $\frac{1}{\rho}$, 
\begin{align}
    \xi_{i,k} \sim \text{Normal}(w_{k}\!\cdot\!x_{i}, \frac{1}{\rho}).
\end{align}
This then defines a forward conditional density (or coupling)
\begin{align}
    p(\xi|w)&\propto e^{-\frac{\rho}{2}\sum_{i,k}(\xi_{i,k}-x_{i}\cdot w_{k})^{2}},
\end{align}
and a joint density for $w,\xi$,
\begin{align}
    p(w, \xi)&= p(w)p(\xi|w)
\end{align}
Via Bayes' rule, this joint density also has expression using the induced marginal on the auxiliary $\xi$ random vector and the reverse conditional density on $w|\xi$,
\begin{align}
    p(w,\xi)&= p(\xi)p(w|\xi).
\end{align} 
As we will show, this choice of forward coupling provides a negative definite correction to the Hessian of the log of $p(w|\xi)$ compared to what we had with $p(w)$, resulting in a log-concave reverse conditional density.

\subsection{Reverse Conditional Density $p(w|\xi)$}
First, we allow for $\xi_{i,k}$ to take arbitrary real values arising from the conditional normal distribution.
\begin{theorem}\label{rev_coup_thm_1}
   With the continuous uniform prior and $\xi_{i,k} \sim \text{Normal}(x_{i}\cdot w_{k}, 1/\rho)$ with $\rho\ge \rho_0$, the reverse conditional density $p(w|\xi)$ is log-concave for the given $\xi$ coupling.
\end{theorem}
\begin{proof}
    The log of the reverse conditional density $p(w|\xi)$ is given by
    \begin{align}
        \log p(w|\xi)=& - \beta \ell(w)+B(\xi)\\
        &-\sum_{i=1}^{N}\sum_{k=1}^{K}\frac{\rho}{2}(\xi_{i,k}-w_{k}\cdot x_{i})^{2},\label{quadratic_term}
    \end{align}
    for some function $B(\xi)$ which does not depend on $w$ and is only required to make the density integrate to 1. The term \eqref{quadratic_term} is a negative quadratic in $w$ which treats each $w_{k}$ as an independent normal random variable. Thus, the additional Hessian contribution will be a $(Kd) \times (Kd)$ negative definite block diagonal matrix with $d \times d$ blocks of the form $\rho\sum_{i=1}^{N}x_{i}x_{i}^{\text{\tiny{T}}}$. Denote the Hessian as $H(w|\xi) = \nabla^{2} \log p(w|\xi)$. For any vector $u \in \mathbb{R}^{Kd}$, with blocks $u_{k} \in \mathbb{R}^{d}$, the quadratic form $u^{\text{\tiny{T}}}H(w|\xi)u$ can be expressed as
    \begin{align}
    &-\beta\sum_{i=1}^{N}\Big(\sum_{k=1}^{K} c_{k}\psi'(w_{k}\cdot x_{i})u_{k}\cdot x_{i}\Big)^{2}\\
    &+\sum_{k=1}^{K}\sum_{i=1}^{N}(u_{k}\cdot x_{i})^{2}[ \beta \text{res}_{i}(w)c_{k}\psi''(w_{k}\cdot x_{i})-\rho].\label{hessian_w_xi}
    \end{align}
    By the assumptions on the second derivative of $\psi$ and the definition of $\rho$ we have
    \begin{align}
        \max_{i,k}(\beta \text{res}_{i}(w)c_{k}\psi''(w_{k}\cdot x_{i})-\rho)&\leq 0.
    \end{align}
    So all the terms in the sum in \eqref{hessian_w_xi} are negative. Thus, the Hessian of the log-likelihood of $p(w|\xi)$ is negative semi-definite and $p(w|\xi)$ is a log-concave density.
\end{proof}

%It is to be noted that this argument for the log concavity of $p(w|\xi)$ works for any log-concave prior, not just the uniform on the $\ell_1$ ball.  Also, it can be made to work for various multimodal likelihoods $e^{-\beta \ell(w)}$, including various nonlinear regressions, not just those that arise in fitting neural networks.  With the centering of the auxiliary variables at the $w_k \cdot x_i$, it applies to models in which the dependence of the log-likelihood $\ell(w)$ on $w$ is either expressed via a smooth function $g(w_{1} \cdot x)$ of one such linear combination , sometimes called a single-index model, as in  generalized linear models in which the link function may produce a non-concave log-likelihoods \cite{nelder1972generalized, mccullagh1989generalized} or expressed via a multi-index regression model $g(w_1 \cdot x,...,w_K \cdot x)$ where-in certain choices of $g$ such as linear combinations of ridge-spline bases, create models that correspond to projection pursuit regression \cite{friedman1981projection, huber1985projection}. The single-hidden-layer network model $\sum_k c_k \psi_k(w_k \cdot x)$ is of this latter form with a fixed $g(z_1,...,z_K)$ being the additive function $\sum_k c_k \psi_k(z_k)$.

While this proof offers a simple way to make a conditional density $p(w|\xi)$ which is log-concave, we also wish to study if there is log-concavity of the induced marginal of $p(\xi)$. The joint log-likelihood for $p(w, \xi)$ contains a bilinear term in $\xi, w$ from expanding the quadratic,
\begin{align}
\sum_{k=1}^{K}\sum_{j=1}^{d}w_{k,j}\sum_{i=1}^{N}\xi_{i,k}x_{i,j}\label{bilinear_term}.
\end{align}
We want some control on how large this term can become, so we restrict the allowed support of $\xi$. When we make this support restriction we will need the parameter $\rho$ to be greater than $\rho_0$.
We define a convenience choice of 
\begin{align}
    \rho&= 2 a_{2}\frac{\beta C_{N} V}{K},\label{larger_rho_def}
\end{align}
which we recognize as $\rho=2 \rho_0$.  

For $0<\alpha \le 1/2$, let $z_\alpha\ge 0$ denote the upper $\alpha$ quantile of the standard normal distribution function. This $z_\alpha$ is the value such that $1-\Phi(z_\alpha) = \alpha$, with the familiar inequality $z_\alpha \le \sqrt{2 \log(1/\alpha)}$. 

For a positive $\delta< 1$, we define a constrained set,
\begin{align}
   B&= \Big\{\xi: \max_{j, k}|\sum_{i=1}^{N}x_{i,j}\xi_{i,k}| \leq N + z_{\delta/2Kd} \sqrt{\frac{N}{\rho}}\Big\}.\label{B_set_def}
\end{align}
where we note that $z_{\delta/2Kd} \le \sqrt{2 \log(2Kd/\delta)}$.
We then define our forward conditional distribution for $p^{*}(\xi|w) = p(\xi|w,B)$ as the normal distribution restricted to the set $B$,
\begin{align}
    p^{*}(\xi|w) =p(\xi|w,B)&=\frac{1_{B}(\xi)p(\xi|w)}{P(\xi \in B|w)} \\
    &=1_{B}(\xi)\frac{\prod_{i=1}^{N}\prod_{k=1}^{K}\big(\frac{\rho}{2\pi} \big)^{\frac{1}{2}}e^{-\frac{\rho}{2}(\xi_{i,k} - x_{i}\cdot w_{k})^{2}}}{\int_{B} \prod_{i=1}^{N}\prod_{k=1}^{K}\big(\frac{\rho}{2\pi}\big)^{\frac{1}{2}}e^{-\frac{\rho}{2}(\xi_{i,k} - x_{i}\cdot w_{k})^{2}}d\xi}.
\end{align}
In accordance with this constrained density, the term \eqref{bilinear_term} will be bounded for any choice of $\xi \in B$ and $w_{k} \in S^{d}_{1}$, which will be a useful property in later proofs.

The denominator of this fraction is the normalizing constant of the density as a result of the restriction to the set $B$. Denote the log normalizing constant as $Z(w) = \log[P(\xi \in B|w)]$,
\begin{align}
    Z(w)&= \log {\int_{B} \prod_{i=1}^{N}\prod_{k=1}^{K}\big(\frac{\rho}{2\pi}\big)^{\frac{1}{2}}e^{-\frac{\rho}{2}(\xi_{i,k} - x_{i}\cdot w_{k})^{2}}d\xi}\label{Z_w_def}.
\end{align}
An equivalent expression for the forward coupling is then
\begin{align}
    p^{*}(\xi|w)&= 1_{B}(\xi)\big(\frac{\rho}{2\pi})^{\frac{NK}{2}}e^{-\frac{\rho}{2}\sum_{i,k}(\xi_{i,k}-x_{i}\cdot w_{k})^{2}-Z(w)}.
\end{align}
This construction also yields for $\xi \in B$ the induced marginal density $p^{*}(\xi)$ with respect to Lebesgue measure, 
\begin{align}
p^{*}(\xi)&= \int p(w)p^{*}(\xi|w)\eta(dw)=1_{B}(\xi) \int p(w)e^{-\frac{\rho}{2}\sum_{i,k}(\xi_{i,k}-x_{i}\cdot w_{k})^{2}-Z(w)}\eta(dw),
\end{align}
and the reverse conditional density $p^{*}(w|\xi)$ with respect to reference measure $\eta$,
\begin{align}
p^{*}(w|\xi)&=\frac{p(w)p^{*}(\xi|w)}{p^{*}(\xi)}= \frac{p(w)e^{-\frac{\rho}{2}\sum_{i,k}(\xi_{i,k}-x_{i}\cdot w_{k})^{2}-Z(w)}}{\int p(w)e^{-\frac{\rho}{2}\sum_{i,k}(\xi_{i,k}-x_{i}\cdot w_{k})^{2}-Z(w)}\eta(dw)}.
\end{align}
Note these densities differ from the $p(w|\xi)$ and $p(\xi)$ defined before without restricting to the set $B$ due to the presence of the $Z(w)$ function.

We can then show two results.  The first is that $Z(w)$ has bounded range, and does not modify the density that much for small $\delta$. Via a Holley-Stroock perturbation argument (Appendix, Section \ref{LSI_app}),
%\ \cite{holley1986logarithmic},
we show that $p^{*}(w|\xi)$ has a reasonably controlled log-Sobolev constant, less than $10K^{2}/(1-\delta)$, even with the original $\rho$. The second result is that the Hessian of $Z(w)$ is small, so that, with the somewhat larger $\rho$ of the present section, we obtain the stronger conclusion that $p^*(w|\xi)$ is log-concave for each $\xi$ in $B$.

This step allows us to construct a \textit{log-concave coupling}, which is to say
\begin{align}
p(w)&= \int p^{*}(w|\xi)p^{*}(\xi)d\xi
\end{align}
where both the conditional $p^*(w|\xi)$ and marginal $p^*(\xi)$ are log-concave.  This conclusion provides a simplified story. Moreover, some results for rapid mixing of MCMC such as Ball Walk and Hit and Run for densities on a constrained set have been stated naturally for densities required to be log-concave, not just log-Sobolev \cite{lovasz2007geometry}.

The restriction of $\xi$ to the set $B$ is the restriction to a very likely set under the unconstrained coupling, in particular we have the following:
\begin{lemma}\label{unconstrained_prob}
    For any weight vector $w$ with $\|w_{k}\|_{1} \leq 1$ the set $B$ in \eqref{B_set_def} has probability using $p(\xi|w)$ satisfying
    \begin{align}
    P(\xi \in B|w) \geq 1-\delta.
    \end{align}
\end{lemma}
\begin{proof}
    See Appendix, Section \ref{app_phi_pfs}.
\end{proof}
Furthermore, the function $Z(w)$ is nearly constant, having small first and second derivative. Therefore, the function has little impact on the log-likelihood.
\begin{lemma}\label{small_derivatives}
For any specified vector $u \in \mathbb{R}^{Kd}$, define the value
\begin{align}
    \tilde{\sigma}^{2}&= \frac{\sum_{i=1}^{N}\sum_{k=1}^{K}(u_{k}\cdot x_{i})^{2}}{\rho}.
\end{align}
For positive values $\delta \leq 0.158$, we then have upper bounds,
    \begin{align}
    |u\cdot  \nabla  Z(w)|&\leq \rho \tilde{\sigma}\frac{\delta}{1-\delta}\big(\sqrt{2\log(1/\delta)}+1\big)
    \end{align}
    and
    \begin{align}
%\    |u^{\text{\tiny{T}}}(\nabla^{2}Z(w))u|&\leq \rho^{2}\tilde{\sigma}^{2}}\big[\frac{\delta}{1-\delta}\Big(2\log(2/\delta)\,-1\Big)\,\,+\frac{\delta^{2}{2\pi} \Big].
    |u^{\text{\tiny{T}}}(\nabla^2 Z(w))u|&\leq \rho^2\tilde{\sigma}^2\frac{\delta}{1-\delta}\Big[2\log(2/\delta)\,+1\,\,+\,\frac{\delta}{1-\delta} (2\log(1/\delta)+3) \Big].
\end{align}
Note both bounds go to $0$ as $\delta \to 0$, and thus can be made arbitrarily small by choice of $\delta$.  
\end{lemma}

\begin{proof}
    See Appendix, Section \ref{app_phi_pfs}.
\end{proof}
The restriction to $\delta\le 0.158$ implies $\delta\le 1-\Phi(1)$ so that the upper quantile $z_\delta$ is at least $1$. More generally, the $+1$ in the gradient bound may be replaced by $1/z_\delta$ and the $+1$ and $+3$ in the Hessian bound may be replaced by $1/z_{\delta/2}$ and $2+1/z_{\delta/2}^2$, respectively. Further refinements are specified in the proof in the appendix. %\In any case these are small contributions to the bound, so we prefer to make the resonable restriction to avoid additional complications to the above expressions.

Thus, with restriction to the set $B$, whose size is determined by $\delta$, and a slightly larger $\rho$, we can give a result that is similar to Theorem \ref{rev_coup_thm_1}. Note this result is for $p^{*}(w|\xi)$ which is distinct from $p(w|\xi)$ due to the presence of the $Z(w)$ function in the log-likelihood and the restriction to $\xi \in B$.

\begin{theorem}\label{rev_coup_thm_2}

    Define the notation
    \begin{align}
%    H_{1}(\delta)&=\frac{2}{\sqrt{2\pi}}\frac{\delta}{1-\delta}\sqrt{2\log \frac{2}{\delta}}\\
       H_{1}(\delta)&= \frac{\delta}{1-\delta}\Big(2\log(2/\delta)\,+1\Big)\\
%    H_{2}(\delta)&=\Big(a_{2}\frac{\beta C_{N}V}{K}\Big)^{2}\frac{1}{2\pi}\frac{\delta^{2}}{(1-\delta)^{2}}.
       H_{2}(\delta)&= \frac{\delta^2}{(1-\delta)^2}\Big(2\log(1/\delta)\,+3\Big)
    \end{align}
and let $H(\delta)$ be defined as $H_1(\delta)+H_2(\delta)$.  For $\rho=2\rho_0$ choose positive $\delta\le 0.054$ that satisfies $H(\delta) \le 1/2$.

Then, for the continuous uniform prior, with $\xi$ restricted to the set $B$ defined by $\delta$, 
%and with $\rho$ as in equation \eqref{larger_rho_def}, 
the reverse conditional density $p^{*}(w|\xi)$ is a log-concave density in $w$, for each $\xi$ in $B$.
\end{theorem}
%{\textbf{Comment: Why the $1/10$ control, larger than the $1/100$ control, even though the $H_2$ tends to be much smaller than the $H_1$?}
\begin{proof}
    See Appendix, Section \ref{app_rev_logcon}.
\end{proof}

%\begin{corollary}\label{delta_cor}
%Positive $\delta \le 0.05$
%will satisfy the condition on $H(\delta)$, 
%producing log-concavity of $p^*_n(w|\xi)$ for each $\xi$ in $B$.
%\end{corollary}

The condition $H(\delta) \le 1/2$ is convenient though not essential. It is tailored to the choice $\rho=2\rho_0$. The proof of Theorem \ref{rev_coup_thm_2} will reveal refined choices of $\rho$, depending on $\delta$, that allows for $H(\delta) < 1$, permitting $0 < \delta < 0.115$, as well as permitting $\rho$ arbitrarily close to $\rho_0$ when $\delta$ is sufficiently small.

By itself, the pairing of a target density $p(w)$ with a suitable normal forward coupling $p^{*}(\xi|w)$ to produce a reverse conditional $p^{*}(w|\xi)$ which is log-concave is not new. As later discussed, the same concept is used in proximal sampling methods and diffusion models. In this work, we go further for our model, obtaining that the induced marginal $p^{*}(\xi)$ is itself log-concave, which we call a log-concave coupling.

\subsection{Marginal Density $p^{*}(\xi)$}

\begin{lemma}
    The score and Hessian of the induced marginal density for $p^{*}(\xi)$ for $\xi \in B$ are expressed as
    \begin{align}
        \partial_{\xi_{i,k}}\log p^{*}(\xi)=&-\rho\,\xi_{i,k}+\rho\,x_{i} \cdot E_{p^{*}}[w_{k}|\xi]\label{score_univariate}\\
        \partial_{ \xi_{i_{1}, k_{1}}, \xi_{i_{2}, k_{2}}} \log p^{*}(\xi)=& -\rho 1\{(i_{1}, k_{1}) = (i_{2}, k_{2})\}\\
        &+\rho^{2}\text{Cov}_{p^{*}}[w_{k_{1}} \cdot x_{i_{1}}, w_{k_{2}} \cdot x_{i_{2}}|\xi].
    \end{align}
    Equivalently in vector form using the $N$ by $d$ data matrix $\textbf{X}$,
    \begin{align}
        \nabla \log p^{*}(\xi)&= \rho\Big( -\xi + E_{p^{*}}\Big[\substack{
            \textbf{X} w_{1}\\
            \ldots\\
            \textbf{X} w_{K}
        }|\xi \Big]\Big) \label{score_exp}\\
    \nabla^{2} \log p^{*}(\xi)&= \rho\Big( -I+\rho\,\text{Cov}_{p^{*}}\Big[\substack{
            \textbf{X} w_{1}\\
            \ldots\\
            \textbf{X} w_{K}
        }|\xi \Big]\Big) \label{hess_exp}.
        \end{align}
\end{lemma}
\begin{proof} 
    The stated results are a consequence of simple calculus, but we will derive them using a statistical interpretation that avoids tedious calculations.
    
 Consider $\xi$ in $B$. The log of the induced marginal for $p^{*}(\xi)$ is equal to the log of the joint density with $w$ integrated out,
    \begin{align}
        \log p^{*}(\xi)= \log \Big(\int p(w)p^{*}(\xi|w)\eta(dw)\Big).
    \end{align}
    Rearranging the Gaussian forward conditional, this can be expressed as a quadratic term in $\xi$ and a term which represents a cumulant generating function plus a constant. Recall $Z(w)$ as defined in equation \eqref{Z_w_def}. Denote the function
    \begin{align}
        h(w)&= -\beta \ell(w)-\frac{\rho}{2}\sum_{i=1}^{N}\sum_{k=1}^{K}(w_{k}\cdot x_{i})^{2} - Z(w),
    \end{align}
    which is the part of the log of the joint density which does not depend on $\xi$. The marginal density can then be expressed for $\xi$ in $B$ as
    \begin{align}
    &\log p^{*}(\xi)=-\frac{\rho}{2}\|\xi\|_{2}^{2}\\
    +&\log \Big(\int p_{0}(w)e^{h(w)}e^{\rho \sum_{i=1}^{N}\sum_{k=1}^{K}\xi_{i,k}w_{k} \cdot x_{i} }\eta(dw)\Big)+C, \label{cumulant_gen_func}
    \end{align}
    for some constant $C$ which makes the density integrate to 1. 
%\ Note that $\xi$ is restricted to have support only on the set $B$, so there is an indicator of the set $B$ we do not write in the expression for simplicity.
    
    It is clear that, with suitable constant adjustment, the term \eqref{cumulant_gen_func} is the cumulant generating function of ($\rho$ times) the random variable $u(w)$ defined by
    \begin{align}
        u(w)&= \xi \cdot \Big(\substack{
                \textbf{X} w_{1}\\
                \ldots\\
                \textbf{X} w_{K}
            }\Big),
    \end{align}
    when $w$ is distributed according to the density proportional to $p_{0}(w)e^{h(w)}$. Thus, the gradient in $\xi$ is the mean of the vector and the second derivative is the covariance, as are standard properties of derivatives of cumulant generating functions, wherein the mean and covariance are taken with respect to the probability density tilted by $e^{\rho u(w)}$, which is recognized to be the reverse conditional $p^{*}(w|\xi)$.
\end{proof}
We highlight two important consequences of this result.
\begin{corollary}\label{score_cor}
The score $\nabla \log p^{*}(\xi)$ is a linear transformation of the vectors of expected values $E_{p^{*}}[w_{k}|\xi]$ from the log-concave reverse conditional $p^{*}(w|\xi)$.
\end{corollary}
\begin{proof}
 This is what is expressed in \eqref{score_univariate} or \eqref{score_exp}.
\end{proof}
%\begin{remark}
%Therefore, while we do not have an explicit closed form expression for 

\noindent \textbf{Remark 1.}
To obtain the score of the marginal density, the $E_{p^{*}}[w_{k}|\xi]$ can be estimated using an MCMC method and thus it is readily available for use. Note that to run an MCMC algorithm such as Metropolis adjusted Langevin diffusion to obtain samples from $p^{*}(\xi)$, its score is needed at each step to update $\xi$.  The conditional means needed for this score at $\xi$ can be computed via  an "inner loop" consisting of an MCMC algorithm to obtain samples from $p^{*}(w|\xi)$ to average. The log-concavity of this reverse conditional provides for the rapid mixing in this inner loop.
%\end{remark}

\begin{corollary}\label{variance_cor}
    The density $p^{*}(\xi)$ is log-concave if for any unit vector $u \in  \mathbb{R}^{NK}$, with blocks $u_{k} \in \mathbb{R}^{N}$, the variance of a particular linear combination of w, namely
    \begin{align}
        v(w)&= \sum_{k=1}^{K}u_{k}^{\text{\tiny{T}}} \textbf{X}w_{k},
    \end{align}
    with respect to the reverse conditional $p^{*}(w|\xi)$ is less than $1/\rho$ ,
    \begin{align}
        \text{Var}_{p^{*}}[v(w)|\xi]\leq 1/\rho,\label{var_cond}
    \end{align}
    for $\xi$ in the convex support set $B$.
\end{corollary}
\begin{proof}
    This is a simple consequence of \eqref{hess_exp}.
\end{proof}
Therefore, to show that $p^{*}(\xi)$ is log-concave we must provide an upper bound on the covariance of $w$ using the reverse conditional density $p^{*}(w|\xi)$. Note that such conditional expectation and conditional covariance representations would also hold using $p(\xi)$, which is defined without conditioning on the set $B$ and thus does not include the $Z(w)$ term in the joint log density. However, the restrictions imposed on maximum inner products by the definition of $B$ will prove useful in upper bounding the reverse conditional covariance.

\subsection{Conditional Covariance Control}
The log of $p^{*}(w|\xi)$ is the log of the prior density plus an additional concave term. Under a log-concave prior, one would expect that adding a concave term to the exponent of an already log-concave density should result in less variance in every direction. Thus one can conjecture the prior covariance would be more than the conditional covariance for any conditioning value,
\begin{align}
    \text{Cov}_{P_{0}}[w] \succ \text{Cov}_{p^{*}}[w|\xi]\quad \forall \xi \in B. \label{cov_claim}
\end{align}
Under a Gaussian prior, such a statement would follow easily from the Brascamp-Lieb inequality \cite[Proposition 2.1]{bobkov2000brunn}. However, for the uniform prior on a convex set, this method does not directly apply.

The covariance matrix of the uniform prior on $(S^{d}_{1})^{K}$ is diagonal (note the different coordinates are uncorrelated but not independent due to symmetry) with entries $\text{Var}_{P_{0}}(w_{k,j}) = \frac{d}{(d+1)^{2}(d+2)} \leq \frac{1}{d^{2}}$ which follows from properties of the Dirichlet distribution. Thus, under conjecture \eqref{cov_claim} we would expect a bound of the form
\begin{align}
    \rho \text{Var}_{p^{*}}[v(w)|\xi]& \leq 2a_{2}\frac{\beta C_{N}V}{K d^{2}}\sum_{j=1}^{d}\sum_{k=1}^{K}(\sum_{i=1}^{N}u_{i,k}x_{i,j})^{2}\\
    & \leq 2a_{2}\frac{\beta C_{N}V N}{K d}\sum_{k=1}^{K}\|u_{k}\|_{1}^{2}\\
    & \leq 2a_{2}\frac{\beta  C_{N}V N}{Kd}\sum_{k=1}^{K}\|u_{k}\|_{2}^{2}\\
    &=  2a_{2}\frac{C_{N}V \beta N}{Kd}.
\end{align}
Thus for $Kd > C(\beta N)$ for some value $C$ we would have log-concavity of the marginal. However, we are unable to prove this conjecture is true. Instead, using a different approach we will conclude for a specified $C$ that
\begin{align}
    Kd\geq C(\beta N)^{2}
\end{align}
results in log-concavity of the marginal density.

Instead of recreating an inequality like \eqref{cov_claim}, we must take a different approach to upper bound the variance in any direction. Denote the function,
    \begin{align}
        h_{\xi}(w)&=-\beta \ell(w)-\sum_{i=1}^{N}\sum_{k=1}^{K}\frac{\rho}{2}(\xi_{i,k}-w_{k}\cdot x_{i})^{2}-Z(w).
    \end{align}
    Denote the function shifted by its prior mean as
    \begin{align}
        \tilde{h}_{\xi}(w)&=h_{\xi}(w)-E_{P_{0}}[h_{\xi}(w)].
    \end{align}
    Define its cumulant generating function with respect to the prior as
    \begin{align}
    \Gamma_{\xi}(\tau)&=\log E_{P_{0}}[e^{\tau \tilde{h}_{\xi}(w)}].\label{Gamma_def}   
    \end{align}

When bounding the variance using the reverse conditional, the expectation can be thought of as taken with respect to the density proportional to the prior times an $e^{\tilde {h}_{\xi}(w)}$ factor. However, because the $h_{\xi}(w)$ is the sum of $N$ varying terms, it is too large to pull out its supremum even for $\xi$ in $B$ and $w$ in $(S_1^d)^K$.  However, with a H{\"o}lder's inequality we can pull out an expectation of $e^{q\tilde {h}_{\xi}(w)}$, which is the expectation of the $q$ power of the factor of interest. Then the variance is controlled when that power $q=\tfrac{\ell}{\ell-1}$  is taken to be close enough to $1$ and when higher order prior moments of $w$ are controlled.  This is captured in the following Lemma and in the subsequent analysis of the two factors that arise from this bound.

    \begin{lemma}\label{holder_main}
        For any integer $\ell \geq 1$ and for any vector $u \in \mathbb{R}^{Kd}$ we have the upper bound
        \begin{align}
            \text{Var}_{p^{*}}(u \cdot w|\xi)& \leq \Big(E_{P_{0}}[(u \cdot w)^{2\ell}]\Big)^{\frac{1}{\ell}}e^{\frac{\ell-1}{\ell}\Gamma_{\xi}(\frac{\ell}{\ell-1})-\Gamma_{\xi}(1)}.\label{holder_eqn}
        \end{align}
    \end{lemma}
    \begin{proof}
        The variance of the inner product $u \cdot w$ is less than its expected square. The reverse conditional density $p^{*}(w|\xi)$ can be expressed as
        \begin{align}
            p^{*}(w|\xi)&= e^{\tilde{h}_{\xi}(w)-\Gamma_{\xi}(1)}p_{0}(w).
        \end{align}
        We then apply a H{\"o}lder's inequality to the integral expression with parameters $p$ and $q$ such that $\frac{1}{p}+\frac{1}{q}=1$,
    \begin{align}
       \text{Var}_{p^{*}}(u \cdot w|\xi)&\leq E_{P_{0}}[(u \cdot w)^{2}e^{\tilde{h}_{\xi}(w)-\Gamma_{\xi}(1)}]\\
       & \leq \Big(E_{P_{0}}[(u \cdot w)^{2p}]\Big)^{\frac{1}{p}}\Big(E_{P_{0}}[e^{q \tilde{h}_{\xi}(w)-q\Gamma_{\xi}(1)}] \Big)^{\frac{1}{q}}.
    \end{align}
    Let $p = \ell$ and $q = \frac{\ell}{\ell-1}$. The second factor can be written as
    \begin{align}
  e^{\frac{\ell-1}{\ell}\Gamma_{\xi}(\frac{\ell}{\ell-1})-\Gamma_{\xi}(1)}.
    \end{align}
\end{proof}

We then study the moments of the prior density and the behavior of the $\Gamma_{\xi}(\tau)$ function separately. As we shall see, when the moments are well controlled, the bound is effective, in part, because it only involves (differences) of cumulants at nearby values.
\begin{lemma}\label{prior_moments_lemma}
    For any unit vector $u \in \mathbb{R}^{NK}$, with blocks $u_{k} \in \mathbb{R}^{N}$, and positive integer $\ell$,
    \begin{align}
        E_{P_{0}}[(\sum_{k=1}^{K}u_{k}^{\text{\tiny{T}}}\textbf{X} w_{k})^{2\ell}]^{\frac{1}{\ell}}&\leq \frac{4\ell N}{\sqrt{e}\,d}.
    \end{align}
\end{lemma}
\begin{proof}
    See Appendix, Section \ref{app_holder_pfs}.
\end{proof}

\begin{lemma}\label{CGF_Growth}
Assume $\rho= 2\rho_0$. Denote the constants
\begin{align}
    A_{1}&=2a_{1}+8a_{2}\label{A_1_def}\\
    A_{2}&=2 \sqrt{2a_{2}}.\label{A_2_def}
\end{align}
Assume positive $\delta \leq  0.116$ and  $d \geq 2$, $K \geq 2$. For any positive integer $\ell \geq 1$ and any $\xi$ from the constrained set $B$, we have
\begin{align}
    &\frac{\ell-1}{\ell}\Gamma_{\xi}(\frac{\ell}{\ell-1}) - \Gamma_{\xi}(1)\leq A_{1}\frac{C_{N}V\beta N}{\ell}+A_{2}\frac{\sqrt{C_{N}V\beta N}}{\ell} \sqrt{2 \log(\frac{2Kd}{\delta})}(\sqrt{1.4K})
\end{align}
\end{lemma}
\begin{proof}
    See Appendix, Section \ref{app_holder_pfs}.
\end{proof}

We summarize the implications of the conclusions of Lemmas 4,5,6 as follows. Ignoring certain constant factors, we have an upper bound on the variance in \eqref{var_cond} for any choice of $\ell$,
\begin{align}
    \frac{N \ell}{d}\,\text{exp}\Big( \frac{\beta N+\sqrt{\beta N K\log(\frac{2Kd}{\delta})}}{\ell} \Big),
\end{align}
which takes the form $(N\ell/d) e^{V_{N,K}/\ell}$, with $V_{N,K}=\beta N+\sqrt{\beta N K\log((2Kd/\delta)}$. Ignoring for now the integer constraint, the optimal continuous choice of $\ell$ to minimize the expression is seen to be $\ell^*=V_{N,K}$, with which the exponent becomes equal to $1$. With this choice, we would have the bound $N V_{N,K} e/d$ on the variance equal to
\begin{align}
    \frac{\beta N^{2} e+N^{\frac{3}{2}}e\sqrt{\beta K{\log(\frac{2Kd}{\delta})}}}{d}.
\end{align}
Multiplying this by $\rho \propto \frac{\beta}{K}$, we would have a bound of order
\begin{align}
    \frac{(\beta N)^{2}}{Kd}\Big(1+\Big[\frac{K\log\Big( \frac{2Kd}{\delta}\Big)}{\beta N}\Big]^{\frac{1}{2}}\Big).
\end{align}
If $K\log(2Kd/\delta) \leq \beta N$, then we have a $O(\frac{(\beta N)^{2}}{Kd})$ bound. With a choice of $d$ and $K$ large enough, we can make this expression be less than 1. We make this statement more precise in the following theorem.

%\ \newpage

\begin{theorem}\label{LogConcaveMarginal}
 Assume $\delta \le 0.116$ and $d\ge 2$, $ K \geq 2$.
%\$d \geq 2, \beta N \geq 2$. 
Further assume that
 \begin{align}
    K\log\Big(\frac{2Kd}{\delta}\Big) \leq \beta N,
 \end{align}
which is essentially a condition that $K$ not be too large, that is, $K$ is less than some multiple of $\beta N/\log(\beta N)$.

With $\rho=2\rho_0$, define $A_{1}, A_{2}$ as in \eqref{A_1_def}, \eqref{A_2_def} and define 
\begin{align}
    A_{3}&=  8\sqrt{e}a_{2}C_{N}V\Big[A_1 C_N V+ 2 A_2 (1.4 C_{N}V)^{\frac{1}{2}}+ \frac{1}{\beta N}\Big].
\end{align}
Let $d$ and $K$ satisfy
\begin{align}
    Kd \geq A_{3}(\beta N)^{2} \label{marginal_log_concave_cond}.
\end{align}
Then $\rho\, \text{Var}_p^*[v(w)|\xi]\le 1$ and the marginal density for $p^{*}(\xi)$ is log-concave with convex support $B$, using the continuous uniform prior. If equation \eqref{marginal_log_concave_cond} is a strict inequality, then the density is strictly log-concave.
\end{theorem}

Log concavity of $p^*(\xi)$ with $\rho=2\rho_0$ is here holding for a larger range of values of $\delta$ than the $\delta \le 0.054$ we require for the log concavity of $p^*(w|\xi)$.  So, with this $\rho$, to get the full log-concave coupling property (with both $p^*(\xi)$ and $p^*(w|\xi)$ log-concave) choose $\delta$ not more than $0.054$.

The combined conditions of Theorem \ref{rev_coup_thm_2} and Theorem \ref{LogConcaveMarginal} are collected in the log-concave coupling Theorem \ref{paper_main_thm} presented earlier in the document.

% A relevant $\delta$ may be $1/Kd$ or a power thereof, though a small constant value such as say $1/300$ is also acceptable (to satisfy Corollary \ref{delta_cor} for example).
\begin{proof}
    By Corollary \ref{variance_cor}, the Hessian of $\log p^{*}(\xi)$ is log-concave when, for any unit vector $u$, we have
    \begin{align}
        \rho \text{Var}_{p^{*}}[\sum_{k=1}^{K}u_{k}^{\text{\tiny{T}}}\textbf{X}w_{k}|\xi] \leq 1.\label{mid_pf_cond}
    \end{align}
    By Lemmas \ref {holder_main}, \ref{prior_moments_lemma}, \ref{CGF_Growth} we have an upper bound for this variance for any scalar $\ell > 1$ and $\xi \in B$. Recall $A_{1}, A_{2}$ as defined in expressions \eqref{A_1_def}, \eqref{A_2_def}. Set $\ell$ to be the least integer greater than or equal to

%    \begin{align}
%        \ell^{*} = A_{1}C_{N}V\beta N+2A_{2}\sqrt{ C_{N}V K\beta N\log(\frac{2Kd}{\delta})}.\label{L_star}
%    \end{align}   
\begin{align}
       \ell^{*} = A_{1}C_{N}V\beta N+2A_{2}\sqrt{ C_{N}V \beta N\log(\frac{2Kd}{\delta})}(\sqrt{1.4K}).\label{L_star}
   \end{align}
This choice makes the cumulant generating difference bound from Lemma \ref{CGF_Growth} not more than $1$, simplifying the variance bound from Lemma \ref{holder_main}.
Using Lemma \ref{prior_moments_lemma}, it gives an upper bound on the product of $\rho$ and variance that appears in \eqref{mid_pf_cond},
    \begin{align}
        &2a_{2}\frac{\beta C_{N} V}{K}\frac{4 N}{\sqrt{e}d}(\ell^{*}+1)\label{apply_L_star}e\\
        &=\, 8\sqrt{ e}A_{1}a_{2}\frac{(C_{N} V \beta N)^{2}}{Kd}\\
        &\,\,+8\sqrt{e}A_{2}a_{2}\frac{2(C_{N} V \beta N)^{\frac{3}{2}}\sqrt{1.4K}}{Kd}\sqrt{\log(\frac{2Kd}{\delta})} \,+\,8\sqrt{ e}a_{2}\frac{C_{N} V \beta N}{Kd}\\
        &\leq\, 8\sqrt{e}a_{2}\frac{(\beta N)^{2}}{Kd} C_N V\Big[
        A_{1}C_{N}V)+2A_2 (1.4C_{N}V)^{\frac{1}{2}}\Big(\frac{K(\log(\frac{2Kd}{\delta})}{\beta N} \Big)^{\frac{1}{2}} + \frac{1}{\beta N} \Big].
    \end{align}
   By assumption,
   \begin{align}
   \frac{K\log(\frac{2Kd}{\delta})}{\beta N} \leq 1,
   \end{align}
   so we have upper bound on \eqref{mid_pf_cond},
   \begin{align}
   8\sqrt{e}a_{2}C_N V \Big[A_1 C_N V+2A_2 (1.4C_{N}V)^{\frac{1}{2}}+ \frac{1}{\beta N}\Big] \frac{(\beta N)^{2}}{Kd}.\label{main_thm_final_exp}
   \end{align}
    
    If $Kd$ satisfies condition \eqref{marginal_log_concave_cond}, then $\rho \text{Var}_{p^*} [v(w)|\xi]$ is less than $1$ for $\xi$ in $B$. By Corollary \ref{variance_cor}, this implies log-concavity of the marginal density $p^*(\xi)$. 
\end{proof}

%\ In the same way the variance $\text{Var}_p[v(w)|\xi]$ can be controlled, but with a smaller constant, when $\xi \ele B$, sufficient to show log-concavity of $p(\xi)$ within $B$.

%\begin{remark}
%    Note that $\ell$ as used in the proof via the H{\"o}lder Inequality must be an integer. This is because we want to consider whole number moments in Lemma \ref{prior_moments_lemma} rather than allowing $\ell$ to be non-whole numbers. Whereas the $\ell^{*}$ in equation \eqref{L_star} is the optimal continuous value. We would have to round up or down to the nearest integer. This would result in $\ell^{*}\pm \epsilon$ for a number $|\epsilon|<1$ in equation \eqref{apply_L_star} instead of $\ell^{*}$. This would give an additional term $\beta N/(Kd)$ in the expression \eqref{main_thm_final_exp}, yet this is a lower order dependence that $(\beta N)^{2}/(Kd)$, so it would still be controlled.
%\end{remark}

%\begin{remark}
\noindent \textbf{Remark 2.}
The log-concavity is obtained in the setting of high interior weight dimension $d$. If the dimension $d$ is not larger than $A_3 (\beta N)^2/K$ then the dimension $d$ can be increased to provide satisfaction of this condition.  We mention two ways to do this. First, the dimension $d$ can be made artificially larger by repeating the input vectors. Say there are original input vectors $\tilde x_{i}$ with default dimension $\tilde{d}$. Define new input vectors by repeating the data $L$ times
\begin{align}
    {x}_{i} &= (\tilde {x}_{i}, \ldots, \tilde{x}_{i}) \in \mathbb{R}^{\tilde{d}L}.
\end{align}
We can then consider ${\textbf{X}}$ as our data matrix with row dimension $d=L \tilde{d}$.

The span of the new data matrix with $\ell_{1}$ controlled weights, $\{z = \textbf{X}w\,:\,\|w\|_{1} \leq 1\}$, is the same as for the original matrix. So we have the same approximation ability of the network. This can also equivalently be considered as inducing a different log-concave prior on the original weight vectors of dimension $\tilde{d}$. Nevertheless, it remains convenient to consider the uniform prior in the higher $d = L\tilde{d}$ dimensional space.  It introduces even more multi-modality in the posterior density $p(w)$ as multiple longer weight vectors yield the same output in the neural network. Yet, with sufficient replication size $L$, by our proceeding theorems the density has a log-concave decomposition. 

As indicated, variable replication does not change what is represented by the network model. A second method to increase the dimensionality incorporates additional model flexibility. Akin to polynomial networks as described in \cite{cover1965geometrical,cover1968capacity,barron1988statistical}, augment the original variables $\tilde x_i$ with a sufficient number of polynomial or spline basis functions of these variables, each bounded by $1$ in absolute value, to create the larger vectors $x_i$ which includes the transformed variables. Computationally excessive dimensionality is avoided by restricting such terms to be additive or first-order interactions.  Retain the activation function with the smoothness conditions we have specified. The network model now represents a larger class of functions of the original variables, yet it retains the desired form as a weighted sum of activation functions applied to linear combinations of (transformed) variables.  The log-concave coupling representation holds for any original dimension $\tilde d$, with sufficiently many of these transformed variables.
%\end{remark}

\par\vspace{0.2cm}
%\begin{remark}
\noindent \textbf{Remark 3.}
This argument for a log-concave coupling works for any log-concave prior on the $\ell_1$ ball whose moments are not more than for the uniform prior. Also, it can be made to work for various likelihoods $e^{-\beta \ell(w)}$ with non-concave $\ell(w)$, including various nonlinear regressions, not just those that arise in fitting neural networks.  In our approach the auxiliary variables are centered at the $w_k \cdot x_i$.  Accordingly, this methodology applies with $K=1$ to models in which the dependence of the log-likelihood $\ell(w)$ on $w$ is expressed via a smooth function $g(w_{1} \cdot x)$ of one such linear combination, sometimes called a single-index model, as in generalized linear models \cite{nelder1972generalized, mccullagh1989generalized} where-in the link function may produce a non-concave log-likelihood. Then these generalized linear models can be provided with a log-concave coupling. With additive transformations of original variables included in the $x$, these represent generalized additive models \cite{hastie1986generalized}, now with a log-concave coupled representation of the posterior distribution.  With $K> 1$ certain multi-index nonlinear regression models $g(w_1 \cdot x,...,w_K \cdot x)$ can be addressed. Particular $g$ such as linear combinations of ridge-spline bases create models that correspond to projection pursuit regression \cite{friedman1981projection, huber1985projection}. Single-hidden-layer network models $\sum_k c_k \psi_k(w_k \cdot x)$ are of this latter form with a fixed $g(z_1,...,z_K)$ being the additive function $\sum_k c_k \psi_k(z_k)$.
%\end{remark}

\section{Discussion}\label{discussion}

The main technique used in this paper which facilitates sampling is the use of a measure decomposition, which is taking a target density $p(w)$ and expressing it in a convenient mixture representation with random variable $\xi$ which may be referred to as a hidden, latent, or auxiliary random variable
\begin{align}
p(w) = \int p(w|\xi)p(\xi)d\xi.
\end{align}
One line of development of this technique has its origins in the statistical physics literature and the study of Ising models or spin glass systems \cite{baker1962certain},\cite{hubbard1972critical}. A tool that  emerges in this literature is the Hubbard transform. For any positive definite matrix $A$, this is simply a rearrangement of a normal density to state
\begin{align}
e^{\frac{1}{2}w^{T}Aw}&= \int \frac{e^{-\frac{1}{2} \xi^{T}A\xi}}{(2\pi)^{d/2} \text{det}(A^{-1})^{1/2}}e^{\xi^{T}A w }d\xi.
\end{align}
Equivalently, the moment generating function of a mean 0 normal random variable is quadratic. As such, whenever a positive definite term of the form $\frac{1}{2}w^{T}Aw$ appears in a log-likelihood, there is an immediate chance to replace it with the Hubbard transform and introduce an auxiliary random variable $\xi$.

This approach has been mostly used as an analysis technique for the mixing time of sampling algorithms for Ising models. Modern works of \cite{bauerschmidt2019very}, \cite{eldan2022spectral} showed significant improvements of the log-Sobolev constant for high temperature spin glass systems with this technique.

Moreover, it has been recognized that a measure decomposition also provides a joint distribution which can be used as a sampling algorithm. For instance, Montanari and Wu \cite{montanari2024provably} use a Gibb’s chain (alternating between conditional $w$ and $\xi$ draws in the present notation) to sample from a Bayesian regression problem with spike and slab prior. This problem has a natural quadratic structure, and thus can make use of the Hubbard transform.

Measure decomposition exists outside of the use of a normal random variable and the Hubbard technique, as studied in more generality in \cite{eldan2018decomposition} and \cite{eldan2020taming}.

More broadly, there are a number of Markov Chain Monte Carlo algorithms that make use of an auxiliary random variable for sampling, without changing the target marginal measure. These include the Swendsen-Wang algorithm with a spin variable and an auxiliary bond variable \cite{swendsen1987nonuniversal}; the Hamiltonian Monte Carlo algorithm with a velocity variable and an auxiliary momentum random variable \cite{betancourt2017conceptual}; the Data Augmentation Gibbs Sampling algorithm arising in Hierarchical Bayes models with latent variables \cite{tanner1987calculation}; the Stochastic Localization algorithm where the Markov chain moves are motivated by a representation of the target random variable through an auxiliary martingale model \cite{el2022sampling, chen2022localization,montanari2023sampling}; and the Proximal Sampling algorithm with a proximal random variable serving as an auxiliary \cite{chen2022improved, huang2024faster}. The Proximal Sampler perhaps bears the most resemblance to our method here. To sample a density $p(w) \propto e^{-V(w)}$, define a new random variable $p(\xi|w) \sim N(w,\frac{1}{\lambda} I)$ for a small value of $\lambda$. Then sample joint density $p(w)p(\xi|w)$. The method is a Gibbs sampler which alternates between sampling the conditional density $p(\xi|w)$ which is a simple normal with small variance, and then the conditional density $p(w|\xi)$ via a rejection algorithm called the Restricted Gaussian Oracle (RGO) \cite{lee2021structured}. Any of these methods may be viewed as a sampling algorithm taking advantage of a specific measure decomposition, or specific joint density.

As we show in this work, the use of measure decomposition is not restricted to densities with a predefined latent structure (for Gibbs sampling) or densities with a quadratic term in them (providing a candidate for the Hubbard transform). Functionally, coupling with a normal forward auxiliary random variable introduces a negative definite corrective term to the Hessian of the log-density, much the same as proximal sampling. Thus we are in essence constructing a measure decomposition ``by force'', defining the forward coupling we want and then deducing the reverse coupling and induced marginal density from the resulting joint density.  The main difference with other techniques is the enforcing of the log-concavity of the reverse coupling and the induced marginal, allowing for provably rapid one pass sampling of each, rather than repeated alternating sampling of the conditionals.

Measure decomposition also has connections with the field of score based diffusion models \cite{song2021scorebased}. Score based diffusions propose starting with a random variable $w'$ from the target density $p(w')$, and then defining the forward SDE $dw_{t}= -w_{t}dt+\sqrt{2}dB_{t}$. At every time $t$, this induces a joint distribution on $p(w', w_{t})$ under which the forward conditional distribution $p(w_{t}|w')$ is a Gaussian distribution with mean being a linear function of $w'$. Paired with this forward SDE is the definition of a reverse SDE that would transport samples from a standard normal distribution to the target distribution of interest. The drift of the reverse diffusion is defined by the scores of the marginal distribution of the forward process $\nabla \log p_{t}(w_{t})$. If these scores can be computed, the target density can be sampled from.

As is the case in our mixture model, the scores of the marginal are defined by expectations with respect to the reverse conditional $p(w'|w_{t})$. For some thresholds $\tau_{1}, \tau_{2}$, for small times $t \leq \tau_{1}$ the reverse conditionals $p(w'|w_{t})$ are log-concave and easily sampled. For large times $t \geq \tau_{2}$, the marginal density $p_{t}(w_{t})$ is approaching a standard normal distribution and thus will become log-concave. If $\tau_{2}< \tau_{1}$, these two regions overlap and the original density $p(w')$ can be written as a log-concave mixture of log-concave components $p(w')= \int p(w'|w_{t})p_{t}(w_{t})dw_{t}$. Thus, the entire procedure of reverse diffusion can be avoided and a one shot sample of $w_{t}$ from its marginal $p_{t}(w_{t})$ and a sample from the reverse conditional $p(w'|w_{t})$ can computed. A variation of this idea is the core procedure we use in this paper, simplifying the processes of a reverse diffusion into one specific and useful choice of joint measure with an auxiliary random variable.

We highlight that the score of the marginal density $\nabla \log p^{*}(\xi)$ is not given as an explicit formula, however it is defined as an expectation with respect to the log-concave reverse conditional for $p^{*}(w|\xi)$ noted in Corollary \ref{score_cor}. Thus, the score of the marginal can be computed as needed via its own MCMC sub-routine.

Here we briefly review sampling literature for log-concave densities. Our density $p^{*}(w|\xi)$ is a weakly log-concave density constrained to a convex set, while $p^{*}(\xi)$ is a strongly log-concave density also restricted to a convex set. For $p^{*}(w|\xi)$, the log of the conditional density only depends on the weight vectors $w$ through their interaction with the data matrix $\textbf{X}w$. The vectors $w$ are $d$ dimensional with $d>N$, thus for any direction orthogonal to the rows of the data matrix the log of the density is flat and has 0 Hessian, hence this density is weakly log-concave. Nonetheless, \cite{lovasz2007geometry} shows Ball Walk and Hit and Run algorithms mix in polynomial time for weakly log-concave densities on a bounded convex set. Recent results in \cite{kook2024sampling} improve upon these mixing time bounds. We note that with a strictly log-concave prior, strict log-concavity of $p^*(w|\xi)$ would be forced. We also note, with different construction of the auxiliary random variable $\xi$, it may be possible to force strict log-concavity in every direction of $p^{*}(w|\xi)$ using a normal with a different mean and covariance matrix for the forward coupling.

In terms of sampling the marginal $p^{*}(\xi)$, we have a strictly log-concave distribution restricted to the convex set defined by $B$. The score $\nabla \log p^{*}(\xi)$ is expressed as a linear transformation of $E^{*}[w|\xi]$ and thus can be computed as needed. If the support set was not restricted, we could use Metropolis Adjusted Langevin Diffusion (MALA) and achieve rapid mixing \cite{dwivedi2019log}. Instead, to deal with the boundary conditions we must use techniques such as a barrier function \cite{srinivasan2024fast} or other adaptations of sampling algorithms to restricted support such as Dikin Walks \cite{kook2024gaussian} and Hamiltonian Monte Carlo in a constrained space \cite{kook2022sampling}.

While in this work we focus on a Bayesian approach and use MCMC for sampling, there have been a number of positive results for training neural networks by optimization in specific instances. For classification problems with well-separated classes and with rather large (potentially overfit) single-hidden-layer networks, \cite{cai2024large} shows that gradient descent with large step size converges quickly to an interpolating solution on the training data (i.e. 0 training loss). The article \cite{tsigler2023benign} demonstrates this solution still has good generalization risk via a form of ``benign overfitting'', however this comes at a cost of being susceptible to adversarial perturbations in specific directions that flip model outputs \cite{frei2024double}.

Another approach to understanding optimization in very large neural networks is to compare them to certain infinite width limits via the Neural Tangent Kernel \cite{jacot2018neural}. With small random initialization and control on the number of gradient steps, gradient methods quickly converge to a near interpolating solution \cite{allen2019convergence,du2019gradient,zou2020gradient}. Networks trained in this regime approach regression with a fixed kernel determined by the covariance of the gradient of the network at the random initialization. This is tantamount to a linear regression onto a prefixed set of basis functions (the large eigenvalue eigenvectors of the kernel).

A related perspective is in \cite{chizat2019lazy} which identifies this regime of an approximate linear model of small weights as ``lazy'' training, that does not allow the internal weights to have much freedom to adjust the basis. Models trained in this regime can have poor generalization, compared to models trained in the more difficult non-lazy regime. In contrast, we seek a procedure that performs as well as if the set of directions $w_{k}$ are adapted to the observed data.

For very wide networks $K > N$, \cite{liu2022loss} shows neural networks satisfy a  Polyak-Łojasiewicz (PL) condition proving convergence of stochastic gradient descent to a global minimizer of the loss function. This is an interesting phenomenon, however without suitable parameter controls (such as $\ell_1$ controls), it is not clear if generalization properties will be favorable in this setting for general function learning.

There are also several negative results \cite{dey2020approximation, goel_et_al:LIPIcs.ITCS.2021.22, froese2024training} showing that training a single-hidden-layer network to interpolation (0 training loss) is an NP hard problem. For example, \cite{vu1998infeasibility} shows that for a network of width $K$, interior weight dimension $d$, and using the step activation function, there does not exist a polynomial time algorithm to achieve average squared training error less than $\zeta (Kd)^{-\frac{3}{2}}$ for an absolute constant $\zeta$. It does not rule out in the noise free setting the possibility of computationally feasible algorithms to achieve average squared error less than a constant times $1/K$.

In this paper, the authors have presented posteriors $p(w)$ that sample all $K$ neuron weights $w_{1}, \ldots, w_{K}$ jointly. However, the problem can also be constructed as a Greedy Bayes procedure sampling one neuron weight at a time based on the residuals of previous fits. The authors discuss these results in \cite{mcdonald2024log_ISIT, barron2024log, barron2024shannon,mcdonald2025computation_thesis}. 

This paper also contains significant portions of Curtis McDonald's PhD thesis, which can be found at \cite{mcdonald2025computation_thesis}. The thesis document contains more in depth analysis of the points discussed here, and a complete look at the greedy approach to sampling.

\section{Conclusion and Future Work}\label{conclusions}

In this work, we study a mixture form of the posterior density for single-hidden-layer neural nets. For a continuous uniform prior on the $\ell_{1}$ ball, we show the posterior density can be expressed as a mixture with only log-concave components when the total number of parameters $Kd$ larger than $C(\beta N)^{2}$ for a constant $C$ where $\beta$ is the inverse temperature and $N$ is the number of data points. 

There are a number of future directions for research. While we are able to give computational control for this problem, no statistical guarantees are derived for this choice of inverse temperature $\beta \leq C \frac{\sqrt{Kd}}{N}$, and a uniform prior. Other choices of prior aside from uniform and couplings different than the normal forward coupling we use here are currently under investigation. These methods may prove better at connecting risk control and a mixture measure form of the target density, but use many of the same ideas outlined in this work.

Further details of the sampling method must also be worked out. The choice of sampling algorithm, hyper-parameter choices such as step size and the number of MCMC iterations, as well as technical details such as condition number have not been addressed in this work. The choice of $\rho$ we make is in a sense the ``smallest'' $\rho$ that forces $p(w|\xi)$ to be log-concave by canceling out any positive definite terms in the Hessian arising from non-linearity (that is, terms dependent on the second derivative of the activation function). Larger choices of $\rho$ can result in stronger log-concavity for the reverse conditional distribution $p(w|\xi)$ that can have sampling benefits.

The H{\"o}lder inequality approach to upper bound the covariance $\text{Cov}[w|\xi]$ is most likely not a tight bound. It is conjectured, for a constant $A$, the covariance of the prior could upper bound the conditional covariance $A \text{Cov}_{P_{0}}[w] \succeq \text{Cov}[w|\xi]$. This would require a lesser condition $Kd> C \beta N$ to achieve log-concavity of $p(\xi)$.

\section{Appendix}

\subsection{\textbf{Proofs for Near Constancy of $Z(w)$}}\label{app_phi_pfs}
In this section, we show the restriction of $\xi$ to the set $B$ is a highly likely event under the base Gaussian distribution, and $Z(w)$ has small magnitude first and second derivatives.

\vspace{0.2cm}
\textbf{Proof of Lemma \ref{unconstrained_prob}}:
\begin{proof}
    We have the event $B=\{\xi: \max_{j,k}|\sum_{i=1}^{N} x_{i,j}\xi_{i,k}| \leq N+z_{\delta/2Kd}\}$. We show that this event is likely for our conditionally independent Gaussian distributions for the $\xi_{i,k}$. This proof follows from standard Gaussian complexity arguments.
    
    The object we must bound is $P(\xi \in B|w)$. If the $\xi_{i,k}$ given $w$ are independent $\text{Normal}(x_{i}\cdot w_{k},1/\rho)$ we may arrange a representation using independent standard normals $Z_{k}$ of dimension $n$,
\begin{align}
    \xi_{k}&= \textbf{X} w_{k}+\frac{1}{\sqrt{\rho}}Z_{k}.
\end{align}

Each mean $x_{i} \cdot w_{k}$ is in $[-1,1]$ due to the weight vector having bounded $\ell_{1}$ norm and the data entries having bounded value.
Consider the complement of the event we want to study. We wish for this event to have probability less than $\delta$.
\begin{align}
    &P\Big[\max_{j,k}|\sum_{i=1}^{N} x_{i,j}\xi_{i,k}| \geq N+z_{\delta/2Kd} \sqrt{\frac{N}{\rho}}\Big],\label{orig_event}
\end{align}
where $P$ is the probability using the normal distribution of $\xi$ given $w$. The max is upper bound by
\begin{align}
    \max_{j,k}|\sum_{i=1}^{N}x_{i,j}\xi_{i,k}| \leq N+\max_{j,k} |\frac{1}{\sqrt{\rho}}\sum_{i=1}^{N}x_{i,j}Z_{i,k}|.
\end{align}
Thus we can bound the larger probability event uniformly for $w \in (S^{d}_{1})^{K}$,
\begin{align}
    P\Big[\max_{j,k}\frac{|\sum_{i=1}^{N} x_{i,j}Z_{i,k}|}{\sqrt{N}} \geq z_{\delta/2Kd}\Big]\leq \delta.
\end{align}
Where the conclusion follows from a union bound and a Gaussian tail bound.
    
\end{proof}
\textbf{Proof of Lemma \ref{small_derivatives}}:
\begin{proof}
    We provide upper bounds on the magnitude of the first and second derivatives of the function $Z(w)$ as defined in equation \eqref{Z_w_def}. Denote $\Phi$ as the standard normal cumulative distribution function (CDF) and $\varphi$ as the standard normal density function. Throughout the proof recall that $p(w|\xi)$ treats each $\xi_{i,k}$ as independent normal with $\xi_{i,k} \sim \text{Normal}(x_{i}\cdot w_{k},\frac{1}{\rho})$ conditionally independent given $w$. The absolute value of the gradient of $Z(w)$ inner product with a vector $u$ with blocks $u_{k}$ is
    \begin{align}
       \big|u\cdot \nabla_{w}Z(w)\big|&= \Big|\rho E[\sum_{i=1}^{N}\sum_{k=1}^{K}(u_{k}\cdot x_{i})(\xi_{i,k} - x_{i}\cdot w_{k})\frac{1_{B}(\xi)}{P(\xi \in B|w)}|w]\Big|.\label{absolute_inner_product}
    \end{align}
%In this representation, the probability of interest is $P(Z \in B_w)$ where $B_w$ is the corresponding set of $Z$ for which $\xi$ is in $B$.
    By Lemma \ref{unconstrained_prob}, the set $B$ has probability at least $1-\delta$. We note the following upper and lower bounds on the Gaussian CDF provided by the classical results of Gordon \cite{gordon1941values},
    \begin{align}
    \frac{\varphi(x)}{x+\frac{1}{x}}\leq 1-\Phi(x) \leq \frac{\varphi(x)}{x}.
    \end{align}
%Consider then the value
%\begin{align}
%\delta^{*}&= \Phi(-\sqrt{2\log(1/\delta)}).
%\end{align}
%For our problem, $Kd \geq 2$ by construction. Then for all positive $\delta\leq 1/e$, it can be shown that $\delta^{*}$ is larger than the term which defines the probability of our set $B$,
%\begin{align}
%\frac{\delta}{\sqrt{2\log(2Kd/\delta)}} \leq \delta^{*}.
%\end{align}
Consider the collections of all measurable sets $D \subset \mathbb{R}^{NK}$ such that $P(\xi \in D) \geq 1-\delta$. This collection contains our original set $B$ as an object in the class. Then, the absolute value of the expected inner product in \eqref{absolute_inner_product} is less than the maximum for any set $D$ in this class,
    \begin{align}
    \max_{\substack{D:\\ P(\xi \in D|w) \geq 1-\delta}} \rho\frac{|E[\sum_{i=1}^{N}\sum_{k=1}^{K}(u_{k}\cdot x_{i})(\xi_{i,k} - x_{i}\cdot w_{k})1_{D}(\xi)|w]|}{1-\delta}.\label{score_term_to_max}
    \end{align}
     Define the value
    \begin{align}
        \tilde{\sigma}&= \sqrt{\frac{\sum_{i=1}^{N}\sum_{k=1}^{K}(u_{k}\cdot x_{i})^{2}}{\rho}}.
    \end{align}
    Under the normal distribution for $\xi$, the integrand in question is a scalar mean 0 normal random variable with this variance,
    \begin{align}
        \sum_{i=1}^{N}\sum_{k=1}^{K}(u_{k}\cdot x_{i})(\xi_{i,k} - x_{i}\cdot w_{k}) \sim \text{Normal}(0,\tilde{\sigma}^{2}).
    \end{align}
Accordingly, since this random variable has mean zero, the expectation in \eqref{score_term_to_max} restricted to $D$ is the same as (minus) the expectation restricted to its complement $D^c$, so the quantity of interest is
   \begin{align}
    \max_{\substack{D:\\ P(\xi \in D^c|w) \leq \delta}} \rho\frac{|E[\sum_{i=1}^{N}\sum_{k=1}^{K}(u_{k}\cdot x_{i})(\xi_{i,k} - x_{i}\cdot w_{k})1_{D^c}(\xi)|w]|}{1-\delta}.
    \end{align}
When $0<\delta \le 1/2$, an extremal set $D^c$ is then seen to correspond to the upper tail of this random variable having the specified probability. That is the maximal value of this expression is achieved with
 %\   The set $D$ which maximizes expression \eqref{score_term_to_max} is then the set which controls the size of this integrand,
    \begin{align}
        D^{c}&= \{\xi: \frac{\sum_{i=1}^{N}\sum_{k=1}^{K}(u_{k}\cdot x_{i})(\xi_{i,k} - x_{i}\cdot w_{k})}{\tilde{\sigma}} \geq z_\delta\},\label{d_star_1}
    \end{align}
where, as before, $z_\delta$ is the upper $\delta$ quantile of the standard normal with $1-\Phi(z_\delta)= \delta$.
(Because of the absolute value and the symmetry of the normal, we can also equally well take the left tail where the standardized normal is less than $- z_\delta$.)
%\The proper choice of $\tau$ is $\sqrt{2\log(1/\delta)}$
From this extremal $D^c$ we then have the upper bound
    \begin{align}
        |u\cdot \nabla_{w}Z(w)|& \leq \frac{\rho \tilde{\sigma}}{1-\delta} \Big|\int_{z_\alpha}^\infty z \varphi(z)dz\Big|=\frac{\rho \tilde{\sigma}}{1-\delta} \varphi(z_\delta),
    \end{align}
using $-z\varphi(z) = \varphi'(z)$ and the fundamental theorem of calculus. By Gordon's inequality $\varphi(z_\delta) \le (z_\delta+1/z_\delta) (1-\Phi(z_\delta))$ which is not more than $\delta (\sqrt{2 \log 1/\delta} + 1)$ for $\delta \le 1-\Phi(1)$.  This yields an upper bound on our expression of interest,
    \begin{align}
%\        |u\cdot \nabla_{w}Z(w)|&\leq \frac{\rho \tilde{\sigma}}{1-\delta}\frac{\delta}{\sqrt{2\pi}}.
 |u\cdot \nabla_{w}Z(w)|&\leq 
\frac{\rho \tilde{\sigma} \delta}{1-\delta} (\sqrt{2 \log (1/\delta)} + 1),
    \end{align}
 which notably goes to 0 as $\delta \to 0$.
    
The Hessian is then a difference in variances,
    \begin{align}
        u^{\text{\tiny{T}}}[\nabla^{2} Z(w)]u=& -\rho \sum_{i=1}^{N}\sum_{k=1}^{K}(x_{i}\cdot u_{k})^{2}\label{bL_t1}\\
        &+\rho^{2}\text{Var}[\sum_{i=1}^{N}\sum_{k=1}^{K}(x_{i}\cdot u_{k})\xi_{i,k}|w,B]\label{bL_t2}.
    \end{align}
Note that the $\xi_{i,k}$ are independent with variance $1/\rho$ conditional on $w$, so if we did not constrain to the set $B$, expressions \eqref{bL_t1} and \eqref{bL_t2} would cancel to $0$. 
%\That is, \eqref{bL_t1} is the variance of the linear function of $\xi$ given $w$ if we did not condition on the set $B$, and \eqref{bL_t2} is the variance conditioned on the set $B$.

The object whose variance we are taking in \eqref{bL_t2} is a linear function of $\xi$, and $\xi$ is a normal random variable given $w$ with diagonal covariance matrix $\frac{1}{\rho}$. By an application of a Brascamp-Lieb inequality, see for example \cite[Proposition 2.1]{bobkov2000brunn}, we would have an upper bound on this variance by the norm square of the vector with entries $x_i \cdot u_k$, divided by $\rho$, which times $\rho^{2}$ is exactly expression \eqref{bL_t1}. Thus, the term \eqref{bL_t2} is less than or equal to the absolute value of term \eqref{bL_t1}. So an upper bound on the quadratic form is $0$, that is $u^{\text{\tiny{T}}}[\nabla^{2}Z(w)]u \leq 0$.

We then compute a lower bound on the variance term in \eqref{bL_t2}. Note a Cramer-Rao lower bound is not applicable here since restriction to a compact set makes integration by parts inapplicable due to boundary conditions. In particular, the expectation of the score of a constrained distribution is not always $0$. 

Using a bias-variance decomposition, write the variance as a non-centered expected squared difference minus the square of the expected difference,
\begin{align}
    &\text{Var}[\sum_{i=1}^{N}\sum_{k=1}^{K}(x_{i}\cdot u_{k})\xi_{i,k}|w,B]\\
    &=  E[\Big(\sum_{i=1}^{N}\sum_{k=1}^{K}(x_{i}\cdot u_{k})(\xi_{i,k}-x_{i}\cdot w_{k})\Big)^{2}|w,\xi \in B]\\
    &-\Big(\sum_{i=1}^{N}\sum_{k=1}^{K}(u_{k}\cdot x_{i})\Big(E[\xi_{i,k}|w,\xi \in B] - x_{i}\cdot w_{k} \Big)\Big)^{2}\label{score_to_be_applied}\\
    & \geq E[\Big(\sum_{i=1}^{N}\sum_{k=1}^{K}(x_{i}\cdot u_{k})(\xi_{i,k}-x_{i}\cdot w_{k})\Big)^{2}\frac{1_{B}(\xi)}{P(\xi \in B|w)}|w]\label{simple_normal}\\
%\    &-\frac{\rho^{2} \tilde{\sigma}^{2}}{(1-\delta)^{2}}\frac{\delta^{2}}{2\pi},\label{score_sq_term}
    &-\frac{ \tilde{\sigma}^{2}\delta^2}{(1-\delta)^{2}}(2 \log (1/\delta) + 3),\label{score_sq_term}
%\NOTE: THE rho^2 is removed here because unlike expression \ref{absolute_inner_product} the quantity we are bounding here does not yet have the rho factor.
\end{align}
where we have applied the previously derived bound on expression \eqref{score_to_be_applied} along with $(z_\delta +1/z_\delta)^2=z_\delta^2 +2+1/z_\delta$, which not more than $z_\delta^2+3$ to deduce expression \eqref{score_sq_term}. 

If we did not condition on the set $B$, the expression \eqref{simple_normal} would be the variance of a simple normal variable with variance $\tilde{\sigma}^{2}$. We will show restricting to $B$ still results in a value very close to $\tilde{\sigma}^{2}$.

The set $B$ has probability at least $1-\delta$.
%The set $B$ has probability at least $1-\delta/\sqrt{2\log(2Kd/\delta)}$.
%Define the value
%\begin{align}
%\delta^{**}&= 2\Phi(-\sqrt{2\log(1/\delta)}).
%\end{align}
%If $Kd \geq 4$, for all positive $\delta\leq 1/16$ we have that $\delta^{**}$ is larger than the term which defines the set $B$ probability,
%\begin{align*}
%\frac{\delta}{\sqrt{2\log(2Kd/\delta)}} \leq \delta^{**}.
%\end{align*}
Then, the expectation in \eqref{simple_normal} restricted to $B$ is lower bounded by the minimum among sets $D$ with $P(\xi \in D) \geq 1-\delta$,
\begin{align}
    &E[\Big(\sum_{i=1}^{N}\sum_{k=1}^{K}(x_{i}\cdot u_{k})(\xi_{i,k}-x_{i}\cdot w_{k})\Big)^{2}\frac{1_{B}(\xi)}{P(\xi \in B|w)}|w]\\
    \geq&\min_{\substack{D:\\P(\xi \in D|w) \geq 1-\delta} } \frac{E[\Big(\sum_{i=1}^{N}\sum_{k=1}^{K}(x_{i}\cdot u_{k})(\xi_{i,k}-x_{i}\cdot w_{k})\Big)^{2}1_{D}(\xi)|w]}{1-\delta}.
\end{align}
The integrand in question, as before, is the same normal variable now squared. The minimizing set $D^{*}$ is then the set placing an upper bound on that expression,
\begin{align}
    D^{*}&=\{\xi:-\tau \leq \frac{\sum_{i=1}^{N}\sum_{k=1}^{K}(x_{i}\cdot u_{k})(\xi_{i,k}-x_{i}\cdot w_{k})}{\tilde{\sigma}} \leq \tau\},
\end{align}
for some value $\tau$, the proper choice being $\tau = z_{\delta/2}$, which achieves the specified probability $1-\delta$.

As an aside, we note this set $D^{*}$ can be deduced from the Neyman-Pearson Lemma \cite[Theorem 3.2.1]{lehmann1986testing}, comparing the distribution where each $\xi_{i,k}$ is independent normal with mean $x_{i} \cdot w_{k}$ and variance $\frac{1}{\rho}$, to the distribution which has this normal density times $(\sum_{i=1}^{N}\sum_{k=1}^{K}(x_{i}\cdot u_{k})(\xi_{i,k}-x_{i}\cdot w_{k}))^{2}$. %\(Likewise, the previous $D^{*}$ in \eqref{d_star_1} can be deduced by a generalization of the Neyman-Pearson lemma in which the alternative is a signed measure measure with the normal density times the factor $\sum_{i=1}^{n}\sum_{k=1}^{K}(x_{i}\cdot u_{k})(\xi_{i,k}-x_{i}\cdot w_{k})$).

We are then integrating a squared normal on a truncated range and have lower bound,
\begin{align}
    &\min_{\substack{D:\\P(\xi \in D|w) \geq 1-\delta} } \frac{E[\Big(\sum_{i=1}^{N}\sum_{k=1}^{K}(x_{i}\cdot u_{k})(\xi_{i,k}-x_{i}\cdot w_{k})\Big)^{2}1_{D}(\xi)|w]}{1-\delta}\\
    &\quad =\quad \frac{\tilde{\sigma}^{2}}{1-\delta} \int_{-z_{\delta/2}}^{z_{\delta/2}}z^{2}\varphi(z)dz.\label{int_express}
\end{align}
To evaluate this integral use its complement set and symmetry of the normal pdf,
\begin{align}
   \int_{-z_{\delta/2}}^{z_{\delta/2}} z^2 \varphi(z)dz&=1-2\int_{z_{\delta/2}}^{\infty} z^2 \varphi(z)dz.
\end{align}
Applying integration by parts, using $(z\varphi(z))'= -z^2 \varphi(z)+ \varphi(z)$, this is
\begin{align}
    1+ 2 z\varphi(z)|_{z_{\delta/2}}^{\infty} -2\int_{z_{\delta/2}}^\infty \varphi(z)dz = 1-2z_{\delta/2}\varphi(z_{\delta/2}) -\delta
\end{align}
Using $z\varphi(z) \le (z^2+1)(1-\Phi(z))$ (from Gordon's inequality), which at $z_{\delta/2}$ is less than $(2\log(2/\delta)+1)\delta/2$, this gives a lower bound for the expression in \eqref{int_express}
\begin{align}
    \frac{\tilde{\sigma}^2}{1-\delta}\Big(1\,-2\delta\big(\log(2/\delta)+1\big)\Big)\label{final_bound_term},
\end{align}
which converges to $\tilde{\sigma}^2$ as $\delta \to 0$. 
We then combine expressions \eqref{bL_t2}, \eqref{score_sq_term}, and \eqref{final_bound_term} to give a lower bound on the Hessian quadratic form,
\begin{align}
%     &u^{\text{\tiny{T}}}[\nabla^{2} Z(w)]u\geq -\rho^{2}\tilde{\sigma}^{2}+\rho^{2}\tilde{\sigma}^{2}(\frac{1}{1-\delta}-\frac{2\delta}{(1-\delta)\sqrt{2\pi}}(\sqrt{2\log(1/\delta)}+\frac{1}{\sqrt{2\log(1/\delta)}}))\\
%     &-\frac{\rho^{2} \tilde{\sigma}^{2}}{(1-\delta)^{2}}\frac{\delta^{2}}{2\pi}\\
%     &= -\frac{\rho^{2}\tilde{\sigma}^{2}}{\sqrt{2\pi}}\frac{\delta}{1-\delta}\Big(-\sqrt{2\pi}+2\sqrt{2\log(1/\delta)}(1+\frac{1}{2\log(1/\delta)})+\frac{1}{\sqrt{2\pi}}\frac{\delta}{1-\delta} \Big)\\
%     &\geq -\frac{\rho^{2}\tilde{\sigma}^{2}}{\sqrt{2\pi}}\frac{\delta}{1-\delta}\Big(2\sqrt{2\log(1/\delta)}+\frac{1}{\sqrt{2\pi}}\frac{\delta}{1-\delta} \Big)
     u^{\text{\tiny{T}}}[\nabla^{2} Z(w)]u&\\
\geq\,
\rho^{2}\tilde{\sigma}^{2}&\Big[-1\,+\,\frac{1}{1-\delta}\big(1- 2\delta(\log(2/\delta)+1)\big)-\frac{\delta^2}{(1-\delta)^2}\big(2\log(1/\delta)+3\big)\Big]\\
= \,
-\rho^{2}\tilde{\sigma}^{2}&\Big[\frac{\delta}{1-\delta}\big(2\log(2/\delta)+1)\big)+\frac{\delta^2}{(1-\delta)^2}\big(2\log(1/\delta)+3\big)\Big]
\end{align}
which converges to 0 as $\delta \to 0$.  This bound completes the proof.

For a refinement, one may use the expressions involving $z_\delta$ and $z_{\delta/2}$ directly.
\begin{align}
     u^{\text{\tiny{T}}}[\nabla^{2} Z(w)]u&\\
\geq\,
\rho^{2}\tilde{\sigma}^{2}\Big[&-1\,+\,\frac{1}{1-\delta}\big(1- 2z_{\delta/2}\varphi(z_{\delta/2})-\delta \big)-\frac{1}{(1-\delta)^2}\big(\varphi(z_\delta)\big)^2\Big]\\
= \,
- \rho^{2}\tilde{\sigma}^2 &\frac{1}{1-\delta} \Big[2z_{\delta/2}\varphi(z_{\delta/2})+\frac{1}{1-\delta}\big(\varphi(z_\delta)\big)^2\Big].
\end{align}
\end{proof}

\subsection{\textbf{Log-Concavity of $p^{*}(w|\xi)$ with Conditioning on $\xi$ in the Set $B$}}\label{app_rev_logcon}
In this section, we show the conditioning of $\xi$ given $w$ to the set $B$ does not affect the log-concavity of the reverse conditional much.

\vspace{0.2cm}
\textbf{Proof of Theorem \ref{rev_coup_thm_2}}
\begin{proof}
We prove the reverse conditional is log-concave when restricting $\xi$ to live in the set $B$. This proof follows much the same way as Theorem \ref{rev_coup_thm_1}.
    The log of $p^{*}(w|\xi)$ is given by
    \begin{align}
        \log p^{*}(w|\xi)=& - \beta \ell(w)+B^*(\xi)\\
        &-\sum_{i=1}^{N}\sum_{k=1}^{K}\frac{\rho}{2}(\xi_{i,k}-w_{k}\cdot x_{i})^{2} \label{quadratic_term_2}\\
       & - Z(w),
    \end{align}
    for some function $B^*(\xi)$ which does not depend on $w$ and is only required to make the density integrate to 1. The term \eqref{quadratic_term_2} is a negative quadratic in $w$ which treats each $w_{k}$ as if it were an independent normal random variable. Thus, the additional Hessian contribution will be a $(Kd) \times (Kd)$ negative definite block diagonal matrix with $d \times d$ blocks of the form $\rho\sum_{i=1}^{N}x_{i}x_{i}^{\text{\tiny{T}}}$. Denote the Hessian as $H(w|\xi)\equiv \nabla^{2} \log p^{*}(w|\xi)$. For any vector $u \in \mathbb{R}^{Kd}$, with blocks $u_{k} \in \mathbb{R}^{d}$, the quadratic form $u^{\text{\tiny{T}}}H(w|\xi)u$ can be expressed as
    \begin{align}
    &-\beta\sum_{i=1}^{N}\Big(\sum_{k=1}^{K} \psi'(w_{k}\cdot x_{i})u_{k}\cdot x_{i}\Big)^{2}\\
    &+\sum_{k=1}^{K}\sum_{i=1}^{N}(u_{k}\cdot x_{i})^{2}\Big[ \beta \text{res}_{i}(w)c_{k}\psi''(w_{k}\cdot x_{i})-\rho)\Big]\label{hessian_w_xi_V2}\\
    &-u^{\text{\tiny{T}}}(\nabla^{2}Z(w))u.\label{bonus_Hessian_term}
    \end{align}
   Let's set $\rho=(1+\epsilon) a_2 C_N V/K$ with positive $\epsilon$, so that it is at least a little bigger than the original $\rho_0$.  With $\epsilon=1$ this is the $\rho=2\rho_0$ assumed in the Theorem statement, but it may be of interest to explore other choices.
    By the assumptions on the second derivative of $\psi$, we have
    \begin{align}
        \max_{i,k}(\beta \text{res}_{i}(w)c_{k}\psi''(w_{k}\cdot x_{i})-\rho)&\leq - \epsilon a_{2}\frac{\beta C_{N}V}{K},
    \end{align}
    so all the terms in the sum in \eqref{hessian_w_xi_V2} are negative. Recall from the definition of $\tilde{\sigma}^{2}$ that
    \begin{align}
    \sum_{k=1}^{K}\sum_{i=1}^{N}(u_{k}\cdot x_{i})^{2} = \rho \tilde{\sigma}^{2}.
    \end{align}
    Therefore, expression \eqref{hessian_w_xi_V2} is less than
    \begin{align}
    -(1+\epsilon)\epsilon \Big(a_{2}\frac{\beta C_{N}V}{K}\Big)^{2}\tilde{\sigma}^{2} = -\frac{\epsilon}{1+\epsilon} \rho^2 \tilde{\sigma}^2.
    \end{align}
    By Lemma \ref{small_derivatives}, for positive $\delta \le 0.158$, a bound on the Hessian term from the correction function $Z$ is
    \begin{align}
        -u^{\text{\tiny{T}}}(\nabla^{2} Z(w))u& \leq \rho^{2}\tilde{\sigma}^{2} \frac{\delta}{1-\delta}\Big(2\log(2/\delta)+1 +\frac{\delta}{1-\delta}\big(2\log(1/\delta)+3\big) \Big).
    \end{align}
We recognize the right side as $\rho^2 \tilde{\sigma}^2 H(\delta)$, where $H(\delta)=H_1(\delta)+H_2(\delta)$ is given by
\begin{align}
H(\delta)=\frac{\delta}{1-\delta}\big(2\log(2/\delta)+1\big) +\frac{\delta^2}{(1-\delta)^2}\big(2\log(1/\delta)+3\big),
\end{align}
an expression that is monotone in $\delta$ for the allowed range and tends to zero as $\delta \to 0$.

    Thus term \eqref{hessian_w_xi_V2} plus \eqref{bonus_Hessian_term} is less than or equal to $\rho^2 \tilde{\sigma}^2$ times the following quantity, which we want to be non-positive,
    \begin{align}
-\frac{\epsilon}{1+\epsilon} + H(\delta)
%        &-\tilde{\sigma}^{2}\Big(a_{2}\frac{\beta C_{N}V}{K}\Big)^{2}\Big(\sqrt{\frac{3}{2}}-1\Big)\Big(\sqrt{\frac{3}{2}}\Big)\label{complex_rev_t1}\\
%        +&\tilde{\sigma}^{2}\Big(a_{2}\frac{\beta C_{N}V}{K}\Big)^{2}\Big(\sqrt{\frac{3}{2}}\Big)^{2}\frac{2}{\sqrt{2\pi}}\frac{\delta}{1-\delta}\sqrt{2\log \frac{\delta}{2}}\\
%        +&\tilde{\sigma}^{2}\Big(a_{2}\frac{\beta C_{N}V}{K}\Big)^{4}\Big(\sqrt{\frac{3}{2}}\Big)^{4}\frac{1}{2\pi}\frac{\delta^{2}}{(1-\delta)^{2}}
\label{complex_rev_t3}
    \end{align}
When $\epsilon=1$ then $\epsilon/(1+\epsilon)$ is $1/2$.  Then by evaluation $-1/2 + H(\delta)$ is seen to be negative for $\delta$ not more than $0.054$, and less than $-1/4$ for $\delta$ not more than $0.0242$. This completes the proof of the log concavity of $p^*(w|\xi)$ for the specified $\rho$ with $\epsilon = 1$.

As long as $H(\delta) < 1$, which is so for positive $\delta < 0.115$, log-concavity can be arranged, by suitable choice of $\epsilon$. Indeed, 
%\optimizing the choice of $\epsilon$, 
we see that $\epsilon = H(\delta)/(1-H(\delta))$ makes the expression $-\frac{\epsilon}{1+\epsilon} + H(\delta)$ equal to zero. With this choice $1+\epsilon = 1/(1-H(\delta))$. This shows $p^*(w|\xi)$ is log-concave with $\rho$ at least $\rho_\delta = \rho_0/(1-H(\delta))$, that is,
\begin{align}
\rho_\delta = \frac{a_2}{1-H(\delta)} \frac{\beta C_N V}{K}
\end{align}
This allows the $\rho$ to be near the original $\rho_0$ for small $\delta$. 

If we use the refined form of the bound on the Hessian of $Z(w)$, the above expression for $H(\delta)$ is replaced by the smaller value
    \begin{align}
H^*(\delta)=  \frac{1}{1-\delta} \Big[2z_{\delta/2}\varphi(z_{\delta/2})+\frac{1}{1-\delta}\big(\varphi(z_\delta)\big)^2\Big].
\end{align}
Accordingly, we can improve the result to obtain log-concavity for a larger range of $\delta$ (as $H^*(\delta)$ remains less than $1$ for $\delta < 0.160$ with no need to impose the restriction that $z_\delta$ be at least $1$), and this log-concavity is obtained for smaller $\rho$, at least $\rho_0/(1-H^*(\delta))$.
\end{proof}

\subsection{\textbf{Proofs of Bounds from the H{\"o}lder Inequality}}\label{app_holder_pfs}
In this section, we bound the two factors in the H{\"o}lder inequality. First, we need a supporting lemma.
\begin{lemma} \label{cont_dir_lemma}
For any vector $x \in [-1,1]^{d}$ and any integer $\ell>0$, the expected inner product with random vector $w$ from the continuous uniform distribution on $S^{d}_{1}$ raised to the power $2\ell$ is upper bound by,
\begin{align}
    E_{P_{0}}[(\sum_{j=1}^{d}x_{j}w_{j})^{2\ell}]&\leq \frac{1}{(d)^{\ell}}\frac{(2\ell)!}{\ell!}.
\end{align}
    
\end{lemma}
\begin{proof}
The sum $\sum_{j=1}^{d}x_{j}w_{j}$ raised to the power $2\ell$ can be expressed as sum using a multi-index $J = (j_{1}, \ldots, j_{2\ell})$ where each $j_{i} \in \{1, \ldots, d\}$ and there are $d^{2\ell}$ terms,
\begin{align}
    E[\Big(\sum_{j=1}^{d}x_{j}w_{j} \Big)^{2\ell}]&= \sum_{j_{1},\ldots,j_{2\ell} }\prod_{i=1}^{2\ell}(x_{j_{i}})E[\prod_{i=1}^{2\ell}w_{j_{i}}].
\end{align}
For a given multi-index vector $J$, let $r(j, J)$ count the number of occurrences of the value $j$ in the vector, $r(j, J)= \sum_{i=1}^{2\ell}1\{j_{i} = j\}$. Then for any multi-index we would have,
\begin{align}
    \prod_{i=1}^{2\ell}w^{j_{i}} = \prod_{j=1}^{d}w_{j}^{r(j, J)}.
\end{align}
Abbreviate $r_{j} = r(j, J)$ for a fixed vector $J$ also note $\sum_{j=1}^{d}r_{j} = 2\ell$. Consider the expectation $E[\prod_{i=1}^{d}w_{j}^{r_{j}}]$. Due to the symmetry of the prior, if any of the $r_{j}$ are odd then the whole expectation is 0. Thus, we only consider vectors $\vec{r} = (r_{1}, \ldots, r_{d})$ where all entries are even.
If we fix the signs of the $w_{j}$ points to live in a given orthant, then the distribution is uniform on the $d+1$ dimensional simplex. Define $w_{d+1} = 1-\sum_{j=1}^{d}|w_{j}|$ then $(|w_{1}|, \ldots, |w_{d}|, w_{d+1})$ has a symmetric Dirichlet $(1, \ldots, 1)$ distribution in $d+1$ dimensions. Note a general Dirichlet distribution in $d+1$ dimensions with parameter vector $\vec{\alpha}=(\alpha_{1}, \ldots, \alpha_{d+1})$ has a properly normalized density as
\begin{align}
    p_{\vec{\alpha}}(w_{1}, \ldots, w_{d})=\frac{\Gamma(\sum_{j=1}^{d}\alpha_{j})}{\prod_{j=1}^{d+1}\Gamma(\alpha_{j})}\prod_{j=1}^{d}(w_{j})^{\alpha_{j}-1}(1-\sum_{j=1}^{d}w_{j})^{\alpha_{d+1}-1}.
\end{align}
Thus the expectation of $\prod_{j=1}^{d}w_{j}^{r_{j}}$ with respect to a symmetric Dirichlet has the form of an un-normalized $\text{Dir}(r_{1}+1, \ldots, r_{d}+1, 1)$ distribution. Thus, the expectation is a ratio of their normalizing constants,
\begin{align}
    E[\prod_{j=1}^{d}w_{j}^{r_{j}}]&= \frac{\Gamma(d+1)\prod_{j=1}^{d}\Gamma(r_{j}+1)}{\Gamma(d+1+\sum_{j=1}^{d}r_{j})} \label{dir_exp}\\
    &= \frac{d!\prod_{j=1}^{d}r_{j}!}{(d+2\ell)!}.
\end{align}
The number of times a specific vector $\vec{r}$ appears from the multi-index $J$ is $\frac{(2\ell)!}{\prod_{j=1}^{d}r_{j}!}$ thus we have,
\begin{align}
    E[(\sum_{j=1}^{d}x_{j}w_{j})^{2\ell}]&= \sum_{\substack{\vec{r}~\text{even} \\\sum_{j}r_{j}=2\ell}}\prod_{j=1}^{d}(x_{j})^{r_{j}}\frac{(2\ell)!}{\prod_{j=1}^{d}r_{j}!}E[\prod_{j=1}^{d}w_{j}^{r_{j}}]\\
    &= \frac{(2\ell!)(d!)}{(d+2\ell)!}\sum_{\substack{\vec{r}~\text{even} \\\sum_{j}r_{j}=2\ell}}\prod_{j=1}^{d}(x_{j})^{r_{j}}\\
    &=\frac{(2\ell!)(d!)}{(d+2\ell)!}\sum_{\substack{\vec{r}~\text{even} \\\sum_{j}r_{j}=2\ell}}\prod_{j=1}^{d}(x_{j}^{2})^{\frac{r_{j}}{2}}\\
    &\leq \frac{(2\ell!)(d!)}{(d+2\ell)!} \frac{(d+\ell-1)!}{\ell!(d-1)!} \label{x_i_bound}\\
    &=\frac{(d+\ell-1)\cdots(d)}{(d+2\ell)\cdots(d+1)}\frac{(2\ell)!}{(\ell)!}\\
    &\leq  \frac{1}{d^{\ell}}\frac{2\ell!}{\ell!},
\end{align}
where inequality \eqref{x_i_bound} follows from each $x_{j}^{2} \leq 1$ thus each term in the sum is less than 1 and there being $\binom{d+\ell-1}{\ell}$ terms in the sum.

\end{proof}

\textbf{Proof of Lemma \ref{prior_moments_lemma}}:
\begin{proof}
We bound the first term in the H{\"o}lder inequality depending on the higher order moments of the prior. We have the unit vector $u \in \mathbb{R}^{nK}$ with n dimensional blocks $u_{k}$. Define vectors in $\mathbb{R}^{d}$ as $v_{k} = \textbf{X}^{\text{\tiny{T}}}u_{k}$ and the object we study is
\begin{align}
    E[(\sum_{k=1}^{K}v_{k} \cdot w_{k})^{2\ell}].
\end{align}
Use a multinomial expansion of this power of a sum and we have expression,
\begin{align}
    E[\sum_{\substack{j_{1}, \ldots, j_{K}\\ \sum j_{k} = 2\ell}}\binom{2\ell}{j_{1},\ldots, j_{K}}\prod_{k=1}^{K}(v_{k}\cdot w_{k})^{j_{k}}]=\sum_{\substack{j_{1}, \ldots, j_{K}\\ \sum j_{k} = 2\ell}}\binom{2\ell}{j_{1},\ldots, j_{K}}\prod_{k=1}^{K}E[(v_{k}\cdot w_{k})^{j_{k}}],
\end{align}
since the prior treats each neuron weigh vector $w_{k}$ as independent and uniform on $S^{d}_{1}$. By the symmetry of the prior, if any $j_{k}$ are odd the whole expression is 0 thus we only sum using even $j_{k}$ values,
\begin{align}
    \sum_{\substack{j_{1}, \ldots, j_{K} \\ \sum j_{k} = \ell}}\binom{2
    \ell}{2j_{1},\ldots, 2j_{K}}\prod_{k=1}^{K}E[(v_{k}\cdot w_{k})^{2j_{k}}]. \label{multi_neuron_sum}
\end{align}
Each vector $v_{k}$ is a linear combination of the rows of the data matrix,
\begin{align}
    v_{k}&= \sum_{i=1}^{N}u_{k,i}x_{i}.
\end{align}
Define $s_{k,i} = \text{sign}(u_{k,i})$ and $\alpha_{k,i} = \frac{|u_{k,i}|}{\|u_{k}\|_{1}}$. We can then interpret the above inner product as a scaled expectation on the data indexes,
\begin{align}
    v_{k} \cdot w_{k}&= (\|u_{k}\|_{1})\sum_{i=1}^{N} \alpha_{k,i} s_{k,i} ~x_{i} \cdot w_{k}.
\end{align}
The average is then less than the maximum term in index $i$,
\begin{align}
    E[(v_{k}\cdot w_{k})^{2j_{k}}]&=(\|u_{k}\|_{1})^{2j_{k}}E[\Big(\sum_{i=1}^{N} \alpha_{k,i} s_{k,i}~x_{i} \cdot w_{k}\Big)^{2j_{k}}]\\
     &\leq (\|u_{k}\|_{1})^{2j_{k}}\sum_{i=1}^{N} \alpha_{k,i}E[\Big(x_{i} \cdot w_{k} \Big)^{2j_{k}}]\\
     &\leq (\|u_{k}\|_{1})^{2j_{k}}\max_{i}E[(x_{i} \cdot w_{k})^{2j_{k}}]\\
     &\leq (\|u_{k}\|_{1})^{2j_{k}}
     \frac{1}{(d)^{j_{k}}}\frac{(2j_{k})!}{j_{k}!},
\end{align}
where we have applied Lemma \ref{cont_dir_lemma}. We then plug this result into equation \eqref{multi_neuron_sum},
\begin{align}
  \frac{1}{d^{\ell}}\frac{(2\ell!)}{\ell!}\Big(\sum_{\substack{j_{1}, \ldots, j_{K} \\ \sum j_{k} = \ell}}
  \binom{\ell}{j_{1},\ldots, j_{K}}\prod_{k=1}^{K}(\|u_{k}\|_{1})^{2j_{k}}\Big)=\frac{1}{d^{\ell}}\frac{(2\ell)!}{\ell!}\Big(\sum_{k=1}^{K}\|u_{k}\|_{1}^{2} \Big)^{\ell}.
\end{align}
For each sub block $u_{k}$ of dimension $n$ we have $\|u_{k}\|_{1}^{2} \leq n \|u_{k}\|^{2}_{2}$ and, since the full vector $u$ is a unit vector, we have $\sum_{k=1}^{K}\|u_{k}\|^{2} = 1$, which gives the upper bound
\begin{align}
    \frac{n^{\ell}(2\ell)!}{d^{\ell} \ell!}.
\end{align}
Via Stirling's bound \cite{robbins1955remark},
\begin{align}
     \sqrt{2\pi \ell} (\frac{\ell}{e})^{\ell}e^{\frac{1}{12\ell+1}}\leq \ell! \leq \sqrt{2\pi \ell} (\frac{\ell}{e})^{\ell}e^{\frac{1}{12\ell}}.
\end{align}
Taking the $\ell$ root we have
\begin{align}
    \Big(\frac{n^{\ell}}{d^{\ell}}\frac{(2\ell)!}{\ell!} \Big)^{\frac{1}{\ell}}& \leq \frac{N}{d} \Big(2^{2\ell+\frac{1}{2}} (\frac{\ell}{e})^{\ell}e^{\frac{1}{24\ell}-\frac{1}{12 \ell+1}}\Big)^{\frac{1}{\ell}}\\
    &= \frac{2^{2+\frac{1}{2\ell}}n \ell}{d} e^{\frac{1}{24 \ell^2}-\frac{1}{12 \ell^2+\ell}-1}\\
    & \leq \frac{4n\ell}{d}\sqrt{2}e^{\frac{1}{24}+\frac{1}{13}-1}\\
    & \leq \frac{4 n\ell}{\sqrt{e} d}.
\end{align}
\end{proof}

\textbf{Proof of Lemma  \ref{CGF_Growth}}:

\begin{proof}
    We bound the second term in the H{\"o}lder inequality determined by the growth rate of the cumulant generating function. By the mean value theorem, there exists some value $\tilde{\tau} \in [1, \frac{\ell}{\ell-1}]$ such that
    \begin{align}
        \Gamma_{\xi}(\frac{\ell}{\ell-1})&= \Gamma_{\xi}(1)+(\Gamma_{\xi})'(\tilde{\tau})[\frac{\ell}{\ell-1}-1].
    \end{align}
Rearranging, we can express the difference
\begin{align}
    \frac{\ell-1}{\ell}\Gamma_{\xi}(\frac{\ell}{\ell-1})-\Gamma_{\xi}(1)=(\Gamma_{\xi})'(\tilde{\tau})\frac{1}{\ell}-\frac{1}{\ell}\Gamma_{\xi}(1).
\end{align}
By construction, $\Gamma_{\xi}(\tau)$ is an increasing convex function with $\Gamma_{\xi}(0)=0$. Thus $\Gamma_{\xi}(1)>0$ and we can study the upper bound
\begin{align}
    \frac{\ell-1}{\ell}\Gamma_{\xi}(\frac{\ell}{\ell-1})-\Gamma_{\xi}(1)\leq(\Gamma_{\xi})'(\tilde{\tau})\frac{1}{\ell}.
\end{align}
Recall $\Gamma_{\xi}(\tau)$ defined in equation \eqref{Gamma_def} is a cumulant generating function of $\tilde{h}^{N}_{\xi}(w)$. Thus, its derivative at $\tilde{\tau}$ is the mean of $\tilde{h}_{\xi}(w)$ under the tilted distribution. The mean is then less than the maximum difference of any two points on the constrained support set,
\begin{align}
    (\Gamma_{\xi})'(\tilde{\tau})&= E_{\tilde{\tau}}[\tilde{h}_{\xi}(w)|\xi]\leq \max_{w, w_{0} \in (S^{d}_{1})^{ K}}(\tilde{h}_{\xi}(w)-\tilde{h}_{\xi}(w_{0})).
\end{align}
By the mean value theorem, for any choice of $w, w_{0}
\in (S^{d}_{1})^{K}$ there exists a $\tilde{w} \in (S^{d}_{1})^{K}$ along the line between $w$ and $w_{0}$ such that
\begin{align}
    \tilde{h}_{\xi}(w)-\tilde{h}_{\xi}(w_{0})&= \nabla_{w}\tilde{h}_{\xi}(\tilde{w})\cdot(w-w_{0}).
\end{align}
For each $k$, the gradient in $w_{k}$, evaluated at $\tilde{w}$ is
\begin{align}
    \nabla_{w_{k}}\tilde{h}_{\xi}(\tilde{w})&= \sum_{i=1}^{N}\big(\beta\text{res}_{i}(\tilde{w})c_{k}\psi'(\tilde {w}_{k}\cdot x_{i}) -\rho w_{k}\cdot x_{i}\big)x_{i}\label{Holder_lemma_2_T1}\\
    &+\rho\beta \sum_{i=1}^{N}\xi_{i,k}x_{i}+\nabla_{w_{k}}Z(w).
\end{align}
Since $|c_k| \le V/K$ and each $|x_i \cdot \tilde{w}_k|\le 1$, with the choice $\rho = 2a_2 \beta C_N V/K$, the multipliers in the sum in \eqref{Holder_lemma_2_T1}, for each $i$, satisfy 
\begin{align}
    \beta \big|\text{res}_{i}(\tilde{w})c_{k}\psi'(w_{k}\cdot x_{i}) -2a_{2}\frac{ C_{N} V}{K}[w_{k}\cdot x_{i}]\big|\leq(a_{1}+2a_{2})\frac{\beta C_{N}V}{K},
\end{align} 
Likewise, we have  $\|x_i \cdot (w_{k}-w_{0,k})\|_{1} \leq 2$. Accordingly, inner product of $w_k -w_{0,k}$ with the first term of the gradient is bounded as
\begin{align}
    \Big[ \beta \sum_{i=1}^{N}(\text{res}_{i}(\tilde{w})c_{k}\psi'(\tilde w_{k}\cdot x_{i}) -2a_2 \frac{C_{N}V}{K})x_{i}\Big]\cdot (w_{k}-w_{0,k})\leq 2\Big(a_{1}+2a_{2}\Big)\frac{C_{N}V \beta N}{K}\label{holder_part2_t1}.
\end{align}
As for the second term, we have 
\begin{align}
    &\rho \Big[\sum_{i=1}^{N}\xi_{i,k}x_{i} \Big]\cdot(w_{k}-w_{0,_{k}}) \leq 2\rho \max_{j} |\sum_{i=1}^{N}\xi_{i,k}x_{i,j}|.
\end{align}
Our original restriction of $\xi$ to the set $B$ is specifically designed to control this term. By definition of the set $B$, for all $k$,
\begin{align}
    \max_{j}|\sum_{i=1}^{N}\xi_{i,k}x_{i,j}| &\leq  N + \sqrt{2 \log (\frac{2Kd}{\delta})}\sqrt{\frac{N}{\rho}}
%\ \\ &= N+\sqrt{2\log \frac{2Kd}{\delta}} \sqrt{\frac{N K}{2a_2\beta C_N V} }.\label{holder_part2_t2}
\end{align}
Multiplying by $2\rho$ it yields the bound on the second term of 
\begin{align}
&2\rho N + 2 \sqrt{\rho N} \sqrt{2 \log \frac{2Kd}{\delta}} \\
&\, = \frac {4a_2 C_N V N}{K} + 2\sqrt{\frac{2a_2 C_N V\beta N}{K}} \sqrt{2 \log \frac{2Kd}{\delta}} \label{holder_part2_t2} 
\end{align}
For the final term, $Z(w)$ is shown to have small derivative. By Lemma \ref{small_derivatives},
\begin{align}
\sum_{k} \nabla_{w_{k}}  Z(w) \cdot (w_{k}-w_{0,k})&\\
\leq \,\, \sqrt{\rho}&\sqrt{\sum_{i=1}^{N}\sum_{k=1}^{K}((w_{k}-w_{0,k})\cdot x_{i})^{2}}\frac{\delta}{1-\delta}\big(2\log(1/\delta)+1\big)\\
\le  2 \sqrt{\rho} & \sqrt{NK} \frac{\delta}{1-\delta}\big(2\log(1/\delta)+1\big)\\
= \,\,2& \sqrt{2a_{2}C_{N}V\beta N}\frac{\delta}{1-\delta}\big(2\log(1/\delta)+1\big) .\label{holder_part2_t3}
\end{align}
Summing using index $k$ for terms \eqref{holder_part2_t1}, \eqref{holder_part2_t2} and combining with term \eqref{holder_part2_t3}, we can upper bound the $(\Gamma_\xi)'(\tilde \tau)$ as,
\begin{align}
%\   &2\Big(a_{1}+2a_{2}\Big)\frac{C_{N}V \beta N}{\ell}+4a_{2}\frac{\beta C_{N}V}{\ell}\Big(N+\sqrt{2\log \frac{2Kd}{\delta}} \sqrt{\frac{N K}{2a_2\beta C_N V}} \Big)\\
%\   &+2\frac{\sqrt{2a_2 C_N V\beta N}}{\ell}\frac{\delta}{1-\delta} \big(2\log(1/\delta) +1\big)\\ =& 
C_{N}V\beta N(2a_1+8a_2)+2\sqrt{2a_2 C_N V\beta N} \Big( \sqrt{2\log \frac{2Kd}{\delta}}\sqrt{K}+    \frac{\delta}{1-\delta}\big(2\log(1/\delta) +1\big)\Big) \label{GammaPrimeBound}
\end{align}
The last part of order $\delta \log 1/\delta$ is negligible compared to the other contributions. With the assumptions $d\ge 2$ and $K \geq 2$,
%, \delta \leq \frac{1}{16}$.  
we find that for all values in the range $0 \le \delta \le 0.116$ we have the inequality
\begin{align}
%\    \frac{\delta}{(1-\delta)} \leq \sqrt{\log \frac{2}{\delta}} \leq \sqrt{\log \frac{2Kd}{\delta}}\sqrt{K}.
%   \frac{\delta}{1-\delta}\big(2\log(1/\delta)+1\big) \leq \sqrt{2\log \frac{2}{\delta}} (2-\sqrt{2}) \leq \sqrt{2 \log \frac{2Kd}{\delta}}(\sqrt{2K}-\sqrt{K}).
 \frac{\delta}{1-\delta}\big(2\log(1/\delta)+1\big) \leq (1/4)\sqrt{2\log \frac{8}{\delta}} \leq (1/4)\sqrt{2 \log \frac{2Kd}{\delta}}
\end{align}
where the form of this bound is chosen for convenience in combining with the second term of (\ref{GammaPrimeBound}). This gives bound on $(\Gamma_\xi)'(\tilde \tau)$ of 
\begin{align}
C_N V\beta N(2a_1+8a_2)+\sqrt{2a_2 C_N V\beta N} \sqrt{2\log \frac{2Kd}{\delta}}(\sqrt{K}+1/4).
\end{align}
Finally, dividing by $\ell$ gives the claimed upper bound on $\frac{\ell-1}{\ell}\Gamma_{\xi}(\frac{\ell}{\ell-1})-\Gamma_{\xi}(1)$
\begin{align}
\frac{1}{\ell} \Big[ C_{N}V\beta N(2a_{1}+8a_{2})+\sqrt{2a_2 C_{N}V\beta N}\Big(\sqrt{\log \frac{2Kd}{\delta}}(\sqrt{K}+1/4\Big)\Big].
\end{align}
When convenient we also use that $\sqrt{K}+1/4$ is not more than $\sqrt{1.4 K}$ for $K\ge 2$.
\end{proof}

If in place of $\rho$ we use the $\rho_\delta = \tfrac{1}{1-H(\delta)} a_2 \tfrac{C_N \beta V}{K}$ defined in the proof of Theorem \ref{rev_coup_thm_2} with the multiplier $\tfrac{1}{1-H(\delta)}$ in place of $2$, then we have a refined bound valid for all $\delta$ with $H(\delta)<1$ (which includes values of $\delta$ less than $0.115$).  This refined bound on $\frac{\ell-1}{\ell}\Gamma_{\xi}(\frac{\ell}{\ell-1})-\Gamma_{\xi}(1)$ is
\begin{align}
\frac{1}{\ell} \Big[ C_{N}V\beta N\big(2a_{1}+4\tfrac{1}{1-H(\delta)} a_{2}\big)+\sqrt{\tfrac{1}{1-H(\delta)} a_2 C_{N}V\beta N}\Big(\sqrt{\log \frac{2Kd}{\delta}}(\sqrt{K}+1/4\Big)\Big],
\end{align}
It provides smaller bounds when $\delta$ is less than $0.054$, that is, when $H(\delta)\le 1/2$. Armed with this $\rho_\delta$, in Lemma \ref{CGF_Growth} and in Theorem \ref{LogConcaveMarginal} we would use 
\begin{align}
A_{1,\delta}  = 2a_{1}+4\frac{1}{1-H(\delta)} a_{2}
\end{align}
and 
\begin{align}
A_{2,\delta} = 2 \sqrt{\frac{1}{1-H(\delta)} a_2}
\end{align}
in place of $A_1$ and $A_2$, respectively. In particular, we would plug these in to provide an improved $A_3$ in Theorem  \ref{LogConcaveMarginal} .

%One can be interested in whether $p(\xi)$ is also log-concave for $\xi$ restricted to $B$. 
While primarily studying the properties of $p^*(\xi)$, it can be of interest to know whether $p(\xi)$ is log-concave at least for the relevant part of its domain when $\xi$ is in $B$ (see the implications below). Accordingly, we can also control the conditional variances using $p(w|\xi)$ instead of $p^*(w|\xi)$.  Toward that end, the analysis above also bounds the second factor of the H{\"o}lder inequality, for $\xi$ in $B$, when the cumulant generating functions uses $p(w|\xi)$.  The only difference is the absence of the $Z(w)$ factor, so that there is not the $|\nabla Z(\tilde w)\cdot (w-w_0)|$ part.  Then there is no  $ \frac{2\delta}{1-\delta}\big(2\log(1/\delta) +1\big)$ part in the above bound, and the final displayed constant is improved with $\sqrt{K}$ in place of $\sqrt{1.4K}$, allowing log-concavity for a slightly larger range of $K$. For certain choices of the parameters it could happen that $p^*(w|\xi)$ or $p^*(\xi)$ are not log-concave while $p(w|\xi)$ and $p(\xi)$ are log-concave for $\xi$ in $B$. When that happens, it is still possible to establish rapid mixing using certain samplers from $p^*(w|\xi)$ and $p^*(\xi)$ via establishment of log-Sobolev properties, as discussed in the next section.

\subsection{Log Sobolev Analysis}\label{LSI_app}

Note that the $Z(w)$ is $\log p(\xi \in B|w)$. Since it is log of a probability it is at most 0. Lemma \ref{unconstrained_prob} gives a lower bound so $Z(w)$ has the following bounded range
\begin{align}
\log \left( 1-\delta\right) \leq Z(w) \leq 0.
\end{align}

Recall the density $p(w|\xi)$ analyzed in Theorem \ref{rev_coup_thm_1} without the $Z(w)$ correction, and using the original $\rho=\rho_0$ (with multiplier $1$). This density $p(w|\xi)$ is weakly log-concave  and it has support contained within the Euclidean ball $B(0,K)$ of radius $K$, since each neuron weight vector has independent bounded 1 norm. As such $p(w|\xi)$ has control on its log-Sobolev constants via standard theory as we here discuss.

Consider all smooth functions $f$, that is all bounded functions with bounded derivatives. For a density $\mu$ the log-Sobolev Constant $C_{LS}(\mu)$ is the smallest constant $C$ such that for all smooth $f$,
\begin{align}
\text{Ent}_{\mu}(f)&= E_{\mu}[f^{2}\ln(f^{2})]- E_{\mu}[f^{2}]\ln(E_{\mu}[f^{2}]) \leq C E_{\mu}[|\nabla f|^{2}] 
\end{align}
Equivalently, the Kullback divergence from $\mu$ of distributions with density w.r.t $\mu$ of the form $g = f^{2}$ is upper bound by a constant times the relative Fisher Information.  The seminal conclusion of Bakry-Emery \cite{bakry1985diffusions,bakry2014analysis} established that strongly log-concave densities have a log-Sobolev constant.  
%\ which they establish while showing, among other properties, an exponential rate of convergence of certain continuous time stochastic diffusion (Langevin diffusion) with a log-concave target.  
Interestingly, when the support is bounded, weakly log-concave densities also have an identified log-Sobolev constant determined by the squared radius of the support.

\begin{lemma}\cite[Corollary 2.3]{bobkov1999isoperimetric}
For any log-concave measure with support contained within the Euclidean ball $B(0,r)$ of radius $r$, the log-Sobolev constant is upper bound by $10 r^{2}$.
\end{lemma}

Control of a log-Sobolev constant is a standard condition for rapid mixing of various MCMC algorithms, so we want to understand the log-Sobolev constant of the measure $p^{*}(w|\xi)$ which includes the $Z(w)$ term. In this section we allow $\rho$ to be the original $\rho_0$, so we are not accounting for the Hessian of $Z(w)$ (so we don't necessarily have log-concavity of $p^*(w|\xi)$). Nevertheless, the $Z(w)$ may be thought of as a small magnitude correction to the otherwise concave exponent, as quantified by Lemma \ref{unconstrained_prob}.  This concave exponent is not necessarily strictly concave, especially when $d> NK$, leading to choices of direction $u$ where the individual $u\cdot x_i$ in the Hessian quadratic form are zero. Nevertheless, the $\ell_1$ constraint sets for the $w_k$ for $k=1,\ldots, K$ lead to $w$ being contained in a Euclidean ball of radius $K$, for which we have the quantified log-Sobolev constant.
% shows for small $\delta$ the range of values $Z(w)$ can take is small. 
As such, a log-Sobolev constant for the distribution with density $p^{*}(w|\xi)$ can be given as an application of \cite{bobkov1999isoperimetric} together with the perturbation lemma of Holley-Stroock \cite[Thm 1.1]{cattiaux2022functional},\cite{holley1987logarithmic}

\begin{lemma}
Since $p^{*}(w|\xi)$ is a bounded perturbation to the original $p(w|\xi)$ density, with ratio between $1-\delta$ and $\tfrac{1}{1-\delta}$ and support in $B(0,K)$, it follows that for each $\xi$ the log-Sobolev constant of this density is upper bound by
\begin{align}
C_{LS}(p^{*}(\cdot |\xi)) \leq \frac{1}{1- \delta} 10 K^{2}.
\end{align}
\end{lemma}

The same strategy can be applied to the analysis of the marginal density $p^*(\xi)$ which has the bounded support set $B$.  Within $B$ this $p^*(\xi)$ can be seen to be a perturbation of the density $p(\xi)$. Indeed in $B$, we have 
\begin{align}
\frac{p^*(\xi)}{p(\xi)}=\frac{ \int p(w) p(\xi|w)e^{-Z(w)}\eta(dw)}{\int p(w) p(\xi|w)\eta(dw)},
\end{align}
which we recognize as $\int p(w|\xi) e^{-Z(w)}\eta(dw)$. Accordingly, this ratio is between $1$ and $\tfrac{1}{1-\delta}$, which shows that perturbation property. As for the log-concavity of $p(\xi)$ on the subset $B$, this holds when $\rho \text{Var}_p [v(w)|\xi] \le 1$. As developed above, this can hold for a slightly larger range of values of $K$ than for the log-concavity of $p^*(\xi)$. 

It follows that when $p(\xi)$ is log-concave on the subset $B$, then $p^*(\xi)$ also has a log-Sobolev constant not more than $10r^2/(1-\delta)$, where $r$ is the radius of the smallest ball containing $B$.
The set $B$ constrains $\xi=(\xi_{i,k}\, : \, i \in [1,N],\, k \in [1,K])$ by requiring each column $\xi_k=(\xi_{1,k},\ldots,\xi_{N,k})'$ to be in the polygonal region determined by $d$ pairs of linear inequality conditions of the form $|\sum_i x_{i,j} \xi_{i,k}|$ being less than a specified value.  There is a pair of such conditions for each $j \in [1,d]$.  So, with $d$ greater than $N$, there is the potential for enough linear constraints that the set $B$ is bounded. However, this depends on the nature of the design $X=(x_{i,j}: i \in [1,N], \, j \in [1,d])$. Rather than tackle that complication for the sake only of obtaining a log-Sobolev property, we have preferred in this paper to be slightly more restrictive regarding the allowed parameters such as $d$ and $K$ and obtain the stronger property of strict log-concavity of $p^*(\xi)$ directly.

\subsection{Alternative Strategy to Show Log-Concavity of $p^*(\xi)$}\label{Alternative_app}

There is a slight variant of our approach to establishing the log-concavity of $p^*(\xi)$ that is worthy to mention. One can show that for $\xi$ in $B$, the variance $\text{Var}_{p^*} [v(w)|\xi]$ is not more than $\text{Var}_p [v(w)|\xi]/(1-\delta)$, using the fact that $p^*(w|\xi)/p(w|\xi) \le 1/(1-\delta)$. 
%together with a variance enlargement from a shift of the conditional mean from what it was with the $p^*$ conditional to what it 
So we get the strong log-concavity we want for $p^*(\xi)$ by a simpler variance calculation showing that $\rho \text{Var}_p [v(w)|\xi] < (1-\delta)$, for $\xi$ in $B$.  This can be established by the assumption $Kd > A_{3,\delta} (\beta N)^2/(1-\delta)$ with a smaller $A_{3,\delta}$ that replaces the $1.4$ factor with $1$. With small $\delta$, the $Kd$ satisfies this property more generally than before. This approach has the advantage that the $Z(w)$ term has been removed for the variance calculation using the $p(w|\xi)$, so that there is no need to control $\nabla Z(w)$ in the variance bound.  Nevertheless, it remains the log-concavity of $p^*(\xi)$ that is being established in this way, and the $Z(w)$ factor is needed in the formulation of the $p^*(w|\xi)$ from which we sample in the log-concave coupling.

\bibliographystyle{emss}
%\bibstyle{emss.bst}
\bibliography{log_concave_journal}
%------

%\begin{thebibliography}{99}

%------ Example for a paper in journal:
% \bibitem{article1}
% A.~Petrunin, Parallel transportation for Alexandrov space with curvature bounded below.
% \emph{Geom. Funct. Anal.} \textbf{8} (1998), no.~1, 123--148
% \Zbl{0903.53045} \MR{1601854}

%------ Example for a book:
% \bibitem{book1}
% W.~P. Ziemer, \emph{Weakly differentiable functions}.
% Grad. Texts in Math. 120,  Springer, New York, 1989
%\Zbl{0692.46022} \MR{1014685}

%------ Example for a paper in a book:
% \bibitem{incollection1}
% J.~S. Milne, Introduction to Shimura varieties.
% In \emph{Harmonic analysis, the trace formula, and Shimura varieties},
% pp. 265--378, Clay Math. Proc. 4,
% American Mathematical Society, Providence, RI, 2005
% \Zbl{1148.14011} \MR{2192012}

%------ Example for a preprint on arXiv:
% \bibitem{preprint1}
% D.~V. Nguyen, S.~K. Chilappagari, M.~W. Marcellin, and B.~Vasic,
% LDPC codes from latin squares free of small trapping sets.
% 2010, \arxiv{1008.4177}

%------ Example for a report:
% \bibitem{report1}
% J.~Schöberl, Commuting quasi-interpolation operators.
% Technical report isc-01-10-math, Texas A\&M University, 2001,
% \url{www.isc.tamu.edu/publications-reports/tr/0110.pdf}

%------ Example for a thesis:
% \bibitem{thesis1}
% E.~Giorgi, \emph{The geometric universe}.
% Ph.D. thesis, University of Maryland, College Park, 2002

%\end{thebibliography}

\end{document}